\input amstex
\documentstyle{amsppt}
%
%
%
\nopagenumbers
\catcode`@=11
\redefine\logo@{\relax}
\catcode`\@=\active
\def\startpage#1{\pageno=#1}
\def\const{\operatorname{const}}
\def\cotan{\operatorname{cotan}}
\def\sign{\operatorname{sign}}
\def\Im{\operatorname{Im}}
\def\Re{\operatorname{Re}}
\def\id{\operatorname{id}}
\def\so{\operatorname{so}}
\pagewidth{12.8cm}
\pageheight{20cm}
\startpage{41}
\Monograph
\topmatter
\title\chapter{3}
Analysis of normality equations in $\Bbb R^2$.
\endtitle
\endtopmatter
\loadbold
\document
\head
\S\,1. Transforming the normality equations to Cartesian coordinates.
\endhead
\rightheadtext{\S\,1. Transforming to Cartesian coordinates.}
     The concept of {\it Newtonian dynamical system admitting the
normal shift of curves in $\Bbb R^2$}, which arose in generalizing
classical Bonnet construction, leads to equations
$$
\aligned
&\alpha_4+\frac{B}{|\bold v|}=0,\\
\vspace{2ex}
&\frac{B\,A}{|\bold v|^2}-\beta_1-\beta_3\,\frac{A}{|\bold v|}-\beta_4\,
\frac{B}{|\bold v|}-\alpha_2+\alpha_3\,\frac{B}{|\bold v|}=0.
\endaligned\hskip -3em
\tag1.1
$$
They were called the {\it weak normality equations}. In two-dimensional
case they are simply called the {\it normality equations}, since they
exhaust complete\footnotemark\ system of normality equations in this
case.\footnotetext{In multidimensional case $n\geqslant 3$ complete
system of normality equations includes weak normality equations and
additional normality equations (see more details in \cite{Shr5}).}
Normality equations \thetag{1.1} are partial differential equations
with respect to components of force field $\bold F$ of Newtonian
dynamical system
\adjustfootnotemark{-1}
$$
\xalignat 2
&\quad\dot\bold r=\bold v,
&&\dot\bold v=\bold F(\bold r,\bold v).
\hskip -2em
\tag1.2
\endxalignat
$$
However, looking at the equations \thetag{1.1}, this is not so obvious.
Force field $\bold F$ of dynamical system \thetag{1.2} are represented
in the equations \thetag{1.1} only through coefficients $A$ and $B$ of
the expansion
$$
\bold F=A\cdot\bold N+B\cdot\bold M,\hskip -2em
\tag1.3
$$
and through parameters $\alpha_2$, $\alpha_3$, $\alpha_4$, $\beta_1$,
$\beta_3$, $\beta_4$, which are defined by the expansions of gradients
$\nabla A$, $\nabla B$, $\tilde\nabla A$, and $\tilde\nabla B$:
$$
\xalignat 2
&\nabla A=\alpha_1\,\bold N + \alpha_2\,\bold M,
&&\tilde\nabla A=\alpha_3\,\bold N + \alpha_4\,\bold M,
\hskip -2em\\
\vspace{-1.2ex}
&&&\tag1.4\\
\vspace{-1.2ex}
&\nabla B=\beta_1\,\bold N + \beta_2\,\bold M,
&&\tilde\nabla B=\beta_3\,\bold N + \beta_4\,\bold M.
\hskip -2em
\endxalignat
$$
Our nearest goal is to transform the equations \thetag{1.1} so that
components of force field $\bold F$ would be represented in explicit
form. From the expansions \thetag{1.3} and \thetag{1.4} we derive
the following expressions for $B$ and $\alpha_4$:
$$
\align
&B=\bigl<\bold F,\,\bold M\bigr>=\sum^2_{i=1}M^i\,F_i,\hskip -2em
\tag1.5\\
\vspace{1ex}
&\alpha_4=\bigl<\tilde\nabla A,\,\bold M\bigr>=\sum^2_{i=1}M^i\,
\tilde\nabla_i A.\hskip -2em
\tag1.6
\endalign
$$
For coefficient $A$ in expansion \thetag{1.3} we have the expression
similar to \thetag{1.5}:
$$
A=\bigl<\bold F,\,\bold N\bigr>=\sum^2_{j=1}N^j\,F_j.\hskip -2em
\tag1.7
$$
Due to the relationships \thetag{1.5}, \thetag{1.6}, and \thetag{1.7}
the first normality equation equation in \thetag{1.1} is rewritten as
$$
\sum^2_{i=1}\left(|\bold v|^{-1}\,F_i+\shave{\sum^2_{j=1}}
\tilde\nabla_i\left(N^j\,F_j\right)\right)M^i=0.\hskip -2em
\tag1.8
$$\par
     Now let's do similar transformations in each term of the second
normality equation in \thetag{1.1}. From \thetag{1.5} and \thetag{1.7}
we derive
$$
\frac{B\,A}{|\bold v|^2}=\sum^2_{i=1}\sum^2_{j=1}|\bold v|^{-2}\,F_i\,
F_j\,N^j\,M^i.\hskip -2em
\tag1.9
$$
For second and fifth terms of the second equation \thetag{1.1} we have
$$
-\beta_1-\alpha_2=-\bigl<\nabla B,\,\bold N\bigr>-\bigl<\nabla A,\,
\bold M\bigr>.
$$
This follows from \thetag{1.4} (see relationships \thetag{2.11} and
\thetag{2.12} in Chapter~\uppercase\expandafter{\romannumeral 2}).
Taking into account \thetag{1.5} and \thetag{1.7} we derive 
$$
-\beta_1-\alpha_2=-\sum^2_{i=1}\sum^2_{j=1}\bigl(\nabla_jF_i+\nabla_iF_j
\bigr)\,M^i\,N^j.\hskip -3em
\tag1.10
$$
Remember that $\nabla_iN^j=0$ and $\nabla_jM^i=0$ (see relationships
\thetag{2.14} in Chapter~\uppercase\expandafter{\romannumeral 2}).
We took into account these relationships in deriving \thetag{1.10}.
\par
    Now let's transform fourth term in the second normality equation
\thetag{1.1}. For the coefficient $\beta_4$ in this term we have
$$
\beta_4=\bigl<\tilde\nabla B,\,\bold M\bigr>=\sum^2_{i=1}\sum^2_{j=1}
\tilde\nabla_j(F_i\,M^i)\,M^j.\hskip -3em
\tag1.11
$$
In expanding brackets under the signs of summation in right hand side
of \thetag{1.11} we take into account the relationships \thetag{2.15}
from Chapter~\uppercase\expandafter{\romannumeral 2}. This yields
$$
-\beta_4\,\frac{B}{|\bold v|}=
\sum^2_{j=1}\sum^2_{r=1}\frac{F_j\,N^j\,F_r\,M^r}{|\bold v|^2}
-\sum^2_{i=1}\sum^2_{j=1}\sum^2_{r=1}\frac{\tilde\nabla_jF_i\,M^i
\,M^j\,F_r\,M^r}{|\bold v|}.
\hskip -3em
\tag1.12
$$
Further calculations are similar to those we have already done.
These calculations do not require special comments:
$$
\gather
-\beta_3\,\frac{A}{|\bold v|}=-\sum^2_{i=1}\sum^2_{j=1}\sum^2_{r=1}
\frac{\tilde\nabla_jF_i\,M^i\,N^j\,F_r\,N^r}{|\bold v|},
\hskip -3em
\tag1.13\\
\vspace{1ex}
\alpha_3\,\frac{B}{|\bold v|}=\sum^2_{i=1}\sum^2_{j=1}\sum^2_{r=1}
\frac{\tilde\nabla_jF_r\,N^r\,N^j\,F_i\,M^i}{|\bold v|}.
\hskip -3em
\tag1.14
\endgather
$$
Let's add the equalities \thetag{1.9}, \thetag{1.10}, \thetag{1.12},
\thetag{1.13}, and \thetag{1.14}. Collecting similar terms, we take
into account the relationship
$$
N^j\,N^r+M^j\,M^r=g^{jr}=\delta^{jr}=\cases 1\text{\ \ for \ }j=r,\\
0\text{\ \ for \ }j\neq r,\endcases
$$
which follows from the fact that vectors $\bold N$ and $\bold M$
form orthonormal frame in standard euclidean metric of the space
$\Bbb R^2$. As a result we have
$$
\gathered
\sum^2_{i=1}\sum^2_{j=1}\left(2\,|\bold v|^{-2}\,F_i\,F_j-\nabla_iF_j
-\nabla_jF_i\right)N^j\,M^i\,+\\
\vspace{1ex}
+\,\sum^2_{i=1}\sum^2_{j=1}\left(\,\shave{\sum^2_{r=1}}\frac{N^r\,N^j\,
\tilde\nabla_jF_r}{|\bold v|}\,F_i-\frac{F^j\,\tilde\nabla_jF_i}
{|\bold v|}\right)M^i=0.
\endgathered\hskip -3em
\tag1.15
$$
Let's denote $v=|\bold v|$, let's change sign in \thetag{1.15} and
join the equations \thetag{1.8} and \thetag{1.15} into one system.
This yields the system of two partial differential equations of the
first order. They are written as follows:
$$
\align
&\sum^2_{i=1}\left(v^{-1}\,F_i+\shave{\sum^2_{j=1}}
\tilde\nabla_i\left(N^j\,F_j\right)\right)M^i=0,\hskip -3em
\tag1.16\\
\vspace{2ex}
&\aligned
&\sum^2_{i=1}\sum^2_{j=1}\left(\nabla_iF_j+\nabla_jF_i-2\,v^{-2}
\,F_i\,F_j\right)N^j\,M^i\,+\\
\vspace{1ex}
&+\sum^2_{i=1}\sum^2_{j=1}\left(\frac{F^j\,\tilde\nabla_jF_i}{v}
-\sum^2_{r=1}\frac{N^r\,N^j\,\tilde\nabla_jF_r}{v}\,F_i\right)M^i=0.
\endaligned\hskip -3em
\tag1.17
\endalign
$$
These are the normality equations for the force field $\bold F(\bold r,
\bold v)$ of Newtonian dynamical system $\ddot\bold r=\bold F(\bold r,
\dot\bold r)$ in $\Bbb R^2$, written in Cartesian coordinates.\par
     In \thetag{1.16} and \thetag{1.17} we have the entries of components
of unitary vector $\bold N$. It is uniquely determined by velocity vector:
$\bold N=|\bold v|^{-1}\cdot\bold v$. Unitary vector $\bold M$ is
orthogonal to $\bold N$, it is also determined by velocity vector (see
formula \thetag{2.7} in Chapter~\uppercase\expandafter{\romannumeral 2}).
However, in order to determine it we should fix some orientation in the
space $\Bbb R^2$. Change of orientation leads to the transformation
$\bold M$ to $-\bold M$. This transformation changes sign of left hand
side of the equations \thetag{1.16} and \thetag{1.17}, but, in essential,
this do not change the equations themselves, since their right hand
side is equal to zero.
    We can make the normality equations \thetag{1.16} and \thetag{1.17}
completely invariant with respect to transformation $\bold M$ to
$-\bold M$. In order to do it we define the operator field $\bold P$ with
the following components:
$$
P^i_j=\delta^i_j-N^i\,N_j=M^i\,M_j.\hskip -3em
\tag1.18
$$
Now, multiplying the equations \thetag{1.16} and \thetag{1.17} by $M_k$,
we write
$$
\align
&\sum^2_{i=1}\left(v^{-1}\,F_i+\shave{\sum^2_{j=1}}
\tilde\nabla_i\left(N^j\,F_j\right)\right)P^i_k=0,\hskip -3em
\tag1.19\\
\vspace{2ex}
&\aligned
&\sum^2_{i=1}\sum^2_{j=1}\left(\nabla_iF_j+\nabla_jF_i-2\,v^{-2}
\,F_i\,F_j\right)N^j\,P^i_k\,+\\
\vspace{1ex}
&+\sum^2_{i=1}\sum^2_{j=1}\left(\frac{F^j\,\tilde\nabla_jF_i}{v}
-\sum^2_{r=1}\frac{N^r\,N^j\,\tilde\nabla_jF_r}{v}\,F_i\right)
P^i_k=0.
\endaligned\hskip -3em
\tag1.20
\endalign
$$
Operator $\bold P$ with components \thetag{1.18} is an operator
of orthogonal projection to the line perpendicular to velocity
vector $\bold v$. It is invariant with respect to the change of
vector $\bold M$ by $-\bold M$.\par
    It's worth to note that the normality equations written as
\thetag{1.19} and \thetag{1.20} coincide by their form with
the weak normality equations for Newtonian dynamical system on
an arbitrary two-dimensional Riemannian manifold $M$ (see
\cite{Bol6} and thesis \cite{Shr5}). The difference of $\Bbb R^2$ and
the case of arbitrary Riemannian manifold $M$ is concentrated in
the definition of differentiations $\nabla$ and $\tilde\nabla$.
This fact will be used below in Chapter~\uppercase
\expandafter{\romannumeral 4}).
\head
\S\,2. Scalar ansatz.
\endhead
    Let's consider first normality equation in the system \thetag{1.1}.
Using this equation, we can express coefficient $B$ in the expansion
\thetag{1.3} through $\alpha_4$:
$$
\pagebreak
B=-|\bold v|\,\alpha_4.\hskip -2em
\tag2.1
$$
Parameter $\alpha_4$ in \thetag{2.1}, in turn, is determined by
coefficient $A$ in the expansion \thetag{1.3} according to the
formula \thetag{1.6}. Therefore we have
$$
B=-|\bold v|\,\sum^2_{i=1}M^i\,\tilde\nabla_i A.\hskip -2em
\tag2.2
$$
Substituting \thetag{2.2} into the expansion \thetag{1.3} we
get formula that expresses components of force vector $\bold F$
through coefficient $A$:
$$
F_k=A\,N_k-|\bold v|\,\sum^2_{i=1}M^i\,M_k\,\tilde\nabla_i A.
\hskip -2em
\tag2.3
$$
If we take into account \thetag{1.18}, then formula \thetag{2.3}
can be transformed to the form
$$
F_k=A\,N_k-|\bold v|\,\sum^2_{i=1}P^i_k\,\tilde\nabla_i A.\hskip -2em
\tag2.4
$$
Formula \thetag{2.4} was called {\bf scalar ansatz}. In \cite{Bol9}
(see also thesis \cite{Shr5}) such formula was used for simplifying
normality equations in multidimensional case $n\geqslant 3$.
Here in two-dimensional case this formula allows us to reduce
system of normality equations \thetag{1.19} and \thetag{1.20} to
one differential equation with respect to the function
$A(\bold r,\bold v)$:
$$
\gathered
\sum^2_{i=1}\left(\nabla_iA+|\bold v|\sum^2_{q=1}\sum^2_{r=1}
P^{qr}\,\tilde\nabla_qA\,\tilde\nabla_r\tilde\nabla_iA\right.-\\
\vspace{0.4ex}
-\left.\shave{\sum^2_{r=1}}N^r\,A\,\tilde\nabla_r\tilde\nabla_iA
-|\bold v|\shave{\sum^2_{r=1}}N^r\,\nabla_r\tilde\nabla_iA\right)
P^i_k=0.
\endgathered\hskip -2em
\tag2.5
$$
Here $P^{qr}=M^q\,M^r$ are components of projection operator $\bold P$.
The equation \thetag{2.5} will be called {\bf reduced normality
equation}. It is as a result of substituting \thetag{2.4} into the
equation \thetag{1.20}. Substituting \thetag{2.4} into the equation
\thetag{1.19} we get the identity. There is no need to reproduce these
calculations here, since they do not differ from those in multidimensional
case (see \cite{Bol9} and thesis \cite{Shr5}).
\head
\S\,3. Polar coordinates in the space of velocities.
\endhead
\parshape 21 0cm 12.8cm 0cm 12.8cm 0cm 12.8cm 0cm 12.8cm
0cm 12.8cm 0cm 12.8cm 0cm 12.8cm 0cm 12.8cm 0cm 12.8cm 0cm 12.8cm 
0cm 12.8cm 4.8cm 8.0cm 4.8cm 8.0cm 4.8cm 8.0cm
4.8cm 8.0cm 4.8cm 8.0cm 4.8cm 8.0cm 4.8cm 8.0cm 4.8cm 8.0cm 4.8cm 8.0cm
0cm 12.8cm
    In reduced normality equation we have the entries of vector $\bold N$
directed along the vector of velocity $\bold v$, and we have the entries
of vector $\bold M$ perpendicular to $\bold v$. This can indicate
the presence of some rotational symmetry in the space of velocities.
In order to use this symmetry for the purposes of further simplifying the
normality equation \thetag{2.5} let's introduce polar coordinates in the
space of velocities. In place of Cartesian coordinates $v^1$ and $v^2$
of velocity vector $\bold v$ we consider radial variable $v=|\bold v|$
and angular variable $\theta$. As $\theta$ we can choose the angle
between velocity vector and some fixed direction in $\Bbb R^2$, for
instance, with the direction positive $OX$ semiaxis. However, we consider
more flexible construction, when $\theta$ is referenced to some unit
vector $\bold m=\bold m(\bold r)$. Let $\tilde\bold m=\tilde
\bold m(\bold r)$ be unit vector perpendicular to $\bold m$:
$$
\xalignat 2
&\bold m=\Vmatrix m^1\\ \vspace{1ex}m^2\endVmatrix,
&&\tilde\bold m=\Vmatrix -m^2\\ \vspace{1ex}m^1\endVmatrix.
\hskip -2em
\tag3.1
\endxalignat
$$
Components of $\bold m$ and $\tilde\bold m$ depend on spatial variables
$r^1$ and $r^2$, \vadjust{\vskip -2.0cm\special{em:graph ris-e21a.gif}
\vskip 2.0cm}but they do not depend on components of $\bold v$. Hence
their velocity gradients are zero:
$$
\xalignat 2
&\tilde\nabla_im^j=0,&&\tilde\nabla_i\tilde m^j=0.\hskip -2em
\tag3.2
\endxalignat
$$
The derivative of any unit vector is perpendicular to this vector.
This yields the following relationships:
$$
\xalignat 2
&\nabla_im^j=\omega_i\,\tilde m^j,
&&\nabla_i\tilde m^j=-\omega_i\,m^j.\hskip -2em
\tag3.3
\endxalignat
$$
Functions $\omega_1(\bold r)$ and $\omega_2(\bold r)$ in \thetag{3.3}
are bound with each other by the relationship 
$$
\nabla_1\omega_2=\frac{\partial\omega^2}{\partial r^1}=\frac{\partial
\omega^1}{\partial r^2}=\nabla_2\omega_1.\hskip -2em
\tag3.4
$$\par
    Let's recalculate the reduced normality equation \thetag{2.5} to
spherical coordinates in velocity space. This means that we should do
the following change of variables in it:
$$
\cases
r^1=x,\\
\vspace{0.2ex}
r^2=y,\\
\vspace{0.2ex}
v^1=v\,(m^1\,\cos\theta+\tilde m^1\,\sin\theta),\\
\vspace{0.2ex}
v^2=v\,(m^2\,\cos\theta+\tilde m^2\,\sin\theta).
\endcases\hskip -2em
\tag3.6
$$
From \thetag{3.6} we obtain the following expressions for vectors
$\bold N$ and $\bold M$:
$$
\aligned
&\bold N=\bold m\cdot\cos\theta+\tilde\bold m\cdot\sin\theta,\\
\vspace{0.4ex}
&\bold M=-\bold m\cdot\sin\theta+\tilde\bold m\cdot\cos\theta.
\endaligned\hskip -2em
\tag3.7
$$
Now let's calculate the derivatives of $\bold v$ in $x$ and $y$
due to \thetag{3.6}:
$$
\pagebreak
\align
&\frac{\partial}{\partial x}\Vmatrix v^1\\ \vspace{1ex} v^2\endVmatrix
=v\,\omega_1\,\cos\theta\,\Vmatrix\tilde m^1\\ \vspace{1ex}\tilde m^2
\endVmatrix-v\,\omega_1\,\sin\theta\,\Vmatrix m^1\\ \vspace{1ex}
m^2\endVmatrix=v\,\omega_1\,\Vmatrix M^1\\ \vspace{1ex} M^2
\endVmatrix,\\
\vspace{1ex}
&\frac{\partial}{\partial y}\Vmatrix v^1\\ \vspace{1ex} v^2\endVmatrix
=v\,\omega_2\,\cos\theta\,\Vmatrix\tilde m^1\\ \vspace{1ex}\tilde m^2
\endVmatrix-v\,\omega_2\,\sin\theta\,\Vmatrix m^1\\ \vspace{1ex}
m^2\endVmatrix=v\,\omega_2\,\Vmatrix M^1\\ \vspace{1ex} M^2
\endVmatrix.
\endalign
$$
These formulas are immediate consequences of \thetag{3.3} and \thetag{3.7}.
The derivatives of $\bold v$ in $v$ and $\theta$ are calculated in a similar
way:
$$
\gather
\frac{\partial}{\partial v}\Vmatrix v^1\\ \vspace{1ex} v^2\endVmatrix
=\cos\theta\,\Vmatrix m^1\\ \vspace{1ex} m^2\endVmatrix+\sin\theta\,
\Vmatrix\tilde m^1\\ \vspace{1ex}\tilde m^2\endVmatrix=\Vmatrix N^1\\
\vspace{1ex} N^2\endVmatrix,\\
\vspace{1ex}
\frac{\partial}{\partial\theta}\Vmatrix v^1\\ \vspace{1ex} v^2\endVmatrix
=v\,\cos\theta\,\Vmatrix\tilde m^1\\ \vspace{1ex}\tilde m^2\endVmatrix-v\,
\sin\theta\,\Vmatrix m^1\\ \vspace{1ex} m^2\endVmatrix=v\,\Vmatrix M^1\\
\vspace{1ex} M^2\endVmatrix.
\endgather
$$
Now we are able to write the Jacoby matrix for the change of variables
\thetag{3.6}:
$$
J=\Vmatrix 1 & 0 & 0 & 0 \\ \vspace{5ex}
0 & 1 & 0 & 0 \\ \vspace{4ex}
v\,\omega_1\,M^1 & v\,\omega_2\,M^1 & N^1 & v\,M^1 \\ \vspace{4ex}
v\,\omega_1\,M^2 & v\,\omega_2\,M^2 & N^2 & v\,M^2
\endVmatrix\hskip -2em
\tag3.8
$$
In order to recalculate differential operators $\nabla_i$ and $\tilde
\nabla_i$ to variables $x$, $y$, $v$, $\theta$ we should find the
matrix inverse to matrix $J$ in \thetag{3.8}:
$$
J^{-1}=\Vmatrix 1 & 0 & 0 & 0 \\ \vspace{4ex}
0 & 1 & 0 & 0 \\ \vspace{4ex}
0 & 0 & N_1 & N_2 \\ \vspace{5ex}
-\omega_1 & -\omega_2 & v^{-1}\,M_1 & v^{-1}\,M_2
\endVmatrix\hskip -2em
\tag3.9
$$
Using formulas \thetag{3.9} for spatial gradients, we obtain
$$
\xalignat 2
&\quad\nabla_1=\frac{\partial}{\partial x}-\omega_1\,\frac{\partial}
{\partial\theta},
&&\nabla_2=\frac{\partial}{\partial y}-\omega_2\,\frac{\partial}
{\partial\theta}.\hskip -3em
\tag3.10
\endxalignat
$$
In a similar way for velocity gradients $\tilde\nabla_1$ and
$\tilde\nabla_2$ we obtain
$$
\pagebreak
\xalignat 2
&\quad\tilde\nabla_1=N_1\,\frac{\partial}{\partial v}+\frac{M_1}{v}\,
\frac{\partial}{\partial\theta},
&&\tilde\nabla_2=N_2\,\frac{\partial}{\partial v}+\frac{M_2}{v}\,
\frac{\partial}{\partial\theta}.\hskip -3em
\tag3.11
\endxalignat
$$
Now we can do the change of variables \thetag{3.6} in the equation
\thetag{2.5}. First of all, let's rewrite this equation itself
$$
\gathered
\sum^2_{i=1}\left(\nabla_iA+|\bold v|\sum^2_{q=1}\sum^2_{r=1}
M^q\,\tilde\nabla_qA\,M^r\,\tilde\nabla_r\tilde\nabla_iA\right.-\\
\vspace{0.4ex}
-\left.\shave{\sum^2_{r=1}}N^r\,A\,\tilde\nabla_r\tilde\nabla_iA
-|\bold v|\shave{\sum^2_{r=1}}N^r\,\nabla_r\tilde\nabla_iA\right)
M^i=0,
\endgathered\hskip -3em
\tag3.12
$$
applying the fact  that $P^i_k=M^i\,M_k$ and $P^{qr}=M^q\,M^r$.
From the relationships \thetag{3.11}, as an immediate consequence,
we have
$$
\xalignat 2
&\quad\sum^2_{r=1}N^r\,\tilde\nabla_r=\frac{\partial}{\partial v},
&&\sum^2_{i=1}M^i\,\tilde\nabla_i=\frac{1}{v}\,\frac{\partial}
{\partial\theta}.\hskip -3em
\tag3.13
\endxalignat
$$
Using the relationships \thetag{3.13} and taking into account formula
\thetag{3.7} for the vector $\bold M$, we calculate third summand
in left hand side of the equation \thetag{3.12}:
$$
-\sum^2_{r=1}\sum^2_{r=1}N^r\,A\,\tilde\nabla_r\tilde\nabla_iA\,M^i=
-A\,\frac{\partial}{\partial v}\!\!\left(\frac{A_\theta}{v}\right)=
-\frac{A\,A_{\theta v}}{v}+\frac{A\,A_\theta}{v^2}.\hskip -3em
\tag3.14
$$
Then consider second summand in left hand side of the equation
\thetag{3.12}:
$$
|\bold v|\sum^2_{q=1}M^q\,\tilde\nabla_qA\sum^2_{i=1}\sum^2_{r=1}\,M^i
\,M^r\,\tilde\nabla_r\tilde\nabla_iA=
\frac{A_\theta}{v}\,\frac{\partial}{\partial\theta}\left(
\shave{\sum^2_{i=1}}M^i\,\tilde\nabla_iA\right)-\frac{A_\theta}{v}
\sum^2_{i=1}\frac{\partial M^i}{\partial\theta}\,\tilde\nabla_iA.
$$
Upon calculating $\partial M^i/\partial\theta$ on the base of \thetag{3.7},
for this summand we get
$$
|\bold v|\sum^2_{q=1}M^q\,\tilde\nabla_qA\sum^2_{i=1}\sum^2_{r=1}\,M^i
\,M^r\,\tilde\nabla_r\tilde\nabla_iA=\frac{A_\theta\,A_{\theta\theta}}{v^2}
+\frac{A_\theta\,A_v}{v}.\hskip -2em
\tag3.15
$$
Now we are to transform first and fourth summands in reduced normality
equation \thetag{3.12}. First of all, note that 
$$
\sum^2_{r=1}N^r\,\nabla_rM^i=0.
$$
This follows from the relationships \thetag{2.14} in
Chapter~\uppercase\expandafter{\romannumeral 2}. Now for the fourth
summand in reduced normality equation \thetag{3.12} we have
$$
\gather
-|\bold v|\sum^2_{i=1}\sum^2_{r=1}M^i\,N^r\,\nabla_r\tilde\nabla_iA=
-v\sum^2_{r=1}N^r\,\nabla_r\left(\frac{A_\theta}{v}\right)=\hskip -3em\\
\vspace{-0.6ex}
\tag3.16\\
\vspace{-1.1ex}
=-N^1\,A_{\theta x}-N^2\,A_{\theta y}+(N^1\,\omega_1+N^2\,\omega_2)\,
A_{\theta\theta}.\hskip -3em
\endgather
$$
In deriving formula \thetag{3.16} the relationships \thetag{3.10} were used.
We use them for to transform first summand in \thetag{3.12} as well:
$$
\sum^2_{i=1}M^i\,\nabla_iA=M^1\,A_x+M^2\,A_y-(M^1\,\omega_1+
M^2\,\omega_2)\,A_\theta.\hskip -1em
\tag3.17
$$
Let's add the relationships \thetag{3.14}, \thetag{3.15}, \thetag{3.16},
\thetag{3.17} and take into account the relationships \thetag{3.7} for
vectors $\bold N$ and $\bold M$ written in coordinates:
$$
\align
\cos\theta&\left(\tilde m^1(A_x-\omega_1\,A_\theta)+\tilde m^2(A_y
-\omega_2\,A_\theta)\,-\right.\hskip -3em\\
&\left.-\,m^1(A_{\theta x}-\omega_1\,A_{\theta\theta})-m^2(A_{\theta y}
-\omega_2\,A_{\theta\theta})\right)\,-\hskip -3em\\
&\hskip -2.1em-\,\sin\theta\left(m^1(A_x-\omega_1\,A_\theta)+m^2(A_y
-\omega_2\,A_\theta)\,+\right.\hskip -3em
\tag3.18\\
&\left.+\,\tilde m^1(A_{\theta x}-\omega_1\,A_{\theta\theta})
+\tilde m^2(A_{\theta y}-\omega_2\,A_{\theta\theta})\right)\,-
\hskip -3em\\
\vspace{2ex}
&\hskip -2.1em-\,\frac{A\,A_{\theta v}}{v}+\frac{A\,A_\theta}{v^2}
+\frac{A_\theta\,A_{\theta\theta}}{v^2}+\frac{A_\theta\,A_v}{v}=0.
\hskip -3em
\endalign
$$
As a result we have transformed the normality equations \thetag{2.5}
to variables $x$, $y$, $v$, $\theta$. Last two variables $v$ and $\theta$
are polar coordinates in the space of velocities. First two variables
$x$ and $y$ are Cartesian coordinates in configuration space $M=\Bbb R^2$.
\head
\S\,4. Some simplest solutions of normality equation.
\endhead
    According to the results of Chapter~\uppercase
\expandafter{\romannumeral 2}, constructing the Newtonian dynamical
systems in $\Bbb R^2$ admitting the normal shift reduces to solving
normality equations \thetag{1.1}. Latter ones are reduced to the
equation \thetag{2.5} for the function $A$. Upon change of variables
\thetag{3.6} the equation \thetag{2.5} turns to \thetag{3.18}.
Formula \thetag{2.2} in variables $x$, $y$, $v$, $\theta$ has the
following very simple form:
$$
B=-\frac{\partial A}{\partial\theta}\hskip -2em
\tag4.1
$$
(see relationships \thetag{3.13} above). For the force field $\bold F$
of Newtonian dynamical system from \thetag{4.1} and from the expansion
\thetag{1.3} we derive
$$
\bold F=A\cdot\bold N-\frac{\partial A}{\partial\theta}\cdot\bold M
\hskip -2em
\tag4.2
$$\par
     Simplest solution of the equation \thetag{3.18} is given by an
identically zero function $A(x,y.v,\theta)=0$. Such solution corresponds
to Newtonian dynamical system with zero force field $\bold F=0$; normal
shift along trajectories of such system coincides with classical
\pagebreak 
Bonnet transformation.\par
     Next (by complexity) solution of the equation \thetag{3.18} is
determined by arbitrary smooth\footnotemark\ function of one variable
$A=A(v)$. \footnotetext{The break of smoothness at the point $v=0$
is admissible, since the equation \thetag{3.18} is written only for
those points, where $v\neq 0$.} If we substitute $A=A(v)$ into
\thetag{4.2}, we obtain $\bold F=A(v)\cdot\bold N$. Force vector
is directed along the vector of velocity, its modulus depends only
on modulus of velocity vector $|\bold v|$:
\adjustfootnotemark{-1}
$$
\bold F=\frac{A(|\bold v|)}{|\bold v|}\cdot\bold v.\hskip -2em
\tag4.3
$$
Force field \thetag{4.3} should also be considered as trivial one.
Trajectories of corresponding dynamical system are straight lines,
normal shift along them coincides with classical Bonnet
transformation.\par
    Less trivial examples of dynamical systems admitting the normal shift
were constructed in paper \cite{Bol3} (see also preprint{Bol2}).
\demo{Example 1} {\bf Spatially homogeneous, but not isotropic force
field}. Let's choose constant vector field $\bold m$ in \thetag{3.1}:
$$
\xalignat 2
&\bold m=\Vmatrix 1\\ \vspace{1ex}0\endVmatrix,
&&\tilde\bold m=\Vmatrix 0\\ \vspace{1ex}1\endVmatrix.
\hskip -2em
\tag4.4
\endxalignat
$$
Then $\omega_1$ and $\omega_2$  in  formulas  \thetag{3.3}  are 
equal
to  zero.  Such  choice  of  $\bold  m$  and  $\tilde\bold   m$ 
substantially
simplifies the normality equation \thetag{3.18}:
$$
\align
(A_y&-A_{\theta x})\,\cos\theta-(A_x+A_{\theta y})\,\sin\theta\,+
\hskip -2em\\
\vspace{-1.0ex}
&\tag4.5\\
\vspace{-1.5ex}
&+\,\frac{A\,A_\theta}{v^2}+\frac{A_\theta\,A_{\theta\theta}}{v^2}
+\frac{A_\theta\,A_v}{v}=\frac{A\,A_{\theta v}}{v}.
\hskip -2em
\endalign
$$
Among all solutions of the equation \thetag{4.5} now we choose those,
which do not depend on $x$ and $y$, i\.~e\. $A=A(v,\theta)$. This
provides spatial homogeneity. Then
$$
\frac{A\,A_\theta}{v^2}+\frac{A_\theta\,A_{\theta\theta}}{v^2}
+\frac{A_\theta\,A_v}{v}=\frac{A\,A_{\theta v}}{v}.\hskip -2em
\tag4.6
$$
Let's take $A=A(v)\,\cos\theta$, where $A(v)$ is an arbitrary smooth
function of one variable. It's easy to check that the function
$A=A(v)\,\cos\theta$ satisfies the equation \thetag{4.6}. Substituting
this function into the formula \thetag{4.2}, we obtain
$$
\bold F=A(v)\cdot(\cos\theta\cdot\bold N+\sin\theta\cdot\bold M).
\hskip -2em
\tag4.7
$$
In Cartesian coordinates force field \thetag{4.7} is written as 
$$
\pagebreak
\bold F=A(|\bold v|)\cdot\left(2\,\bigl<\bold N,\,\bold m\bigr>
\cdot\bold N-\bold m\right)\!,\hskip -2em
\tag4.8
$$
where $\bold m$ is constant vector, which is unitary in our example
and directed along $OX$ axis (see formulas \thetag{4.4}).\par
\parshape 12 4.7cm 8.1cm 4.7cm 8.1cm 4.7cm 8.1cm 4.7cm 8.1cm 4.7cm
8.1cm 4.7cm 8.1cm 4.7cm 8.1cm 4.7cm 8.1cm 4.7cm 8.1cm 4.7cm 8.1cm
4.7cm 8.1cm 0cm 12.8cm
As we see on Fig.~4.1, force vector in force field \thetag{4.7}
form the angle $\theta$ with velocity vector $\bold v$ and it
forms the angle $2\,\theta$ with $OX$ axis. \vadjust{\vskip -0.6cm
\special{em:graph ris-e41a.gif}\vskip 0.6cm}Newtonian dynamical system
with force field \thetag{4.7} admits the normal shift of curves in
the space $\Bbb R^2$. Moreover, it is integrable in quadratures.
In order to show this let's find its trajectories. Force $\bold F$ in
\thetag{4.7} does not depend on coordinates $r^1$ and $r^2$. Therefore
the equation $\ddot\bold r=\bold F(\bold r,\dot\bold r)$ can be written
as $\dot\bold v=\bold F(\bold v)$. Let's integrate this equation in
variables $v$ and $\theta$. Using the relationships \thetag{4.4}, we
write 
$$
\xalignat 2
&\qquad\bold v=v\cdot\Vmatrix \cos\theta\\ \vspace{1ex}\sin\theta\endVmatrix,
&&\dot\bold v=\dot v\cdot\Vmatrix\cos\theta\\ \vspace{1ex}
\sin\theta\endVmatrix+v\,\dot\theta\cdot\Vmatrix -\sin\theta\\
\vspace{1ex}\cos\theta\endVmatrix.
\hskip -2em
\tag4.9
\endxalignat
$$
As for the formula \thetag{4.7} for force vector, we rewrite this
formula as follows:
$$
\bold F=A(v)\,\cos\theta\cdot\Vmatrix\cos\theta\\ \vspace{1ex}
\sin\theta\endVmatrix+A(v)\,\sin\theta\cdot\Vmatrix -\sin\theta\\
\vspace{1ex}\cos\theta\endVmatrix.\hskip -3em
\tag4.10
$$
Comparing \thetag{4.9} and \thetag{4.10}, we see that vectorial 
equation
$\dot\bold v=\bold F(\bold v)$ is written as a  system  of  two 
scalar
equations
$$
\cases \dot v=A(v)\,\cos\theta,\\ \vspace{1ex}v\,\dot\theta=A(v)\,
\sin\theta.\endcases\hskip -3em
\tag4.11
$$
Solution of the system of differential equations \thetag{4.11} is
the pair of functions $v(t)$ and $\theta(t)$. If $A(v)\neq 0$ and
$\sin\theta\neq 0$, then function $\theta=\theta(t)$ is invertible.
The inverse function is $t=t(\theta)$. Substituting $t=t(\theta)$
into $v(t)$, we get the function $v=v(t(\theta))=v(\theta)$. From
the system of two equations \thetag{4.11} we get the following
differential equation for such function:
$$
\frac{dv}{d\theta}=\frac{\dot v}{\vphantom{V^{V^V}}\dot\theta}=
\frac{A(v)\,\cos\theta}{v^{-1}\,A(v)\,\sin\theta}=v\,\cotan\theta.
\hskip -3em
\tag4.12
$$
Differential equation \thetag{4.12} is easily integrable:
$$
v(\theta)=C\,|\sin\theta|\text{, \ where \ }C=\const>0.
\hskip -3em
\tag4.13
$$
Upon substituting \thetag{4.13} into the second equation \thetag{4.11}
we obtain differential equation of the first order for the function
$\theta=\theta(t)$:
$$
\dot\theta=\frac{ A(C\,|\sin\theta|)}{C\,|\sin\theta|}\,\sin\theta.
\hskip -3em
\tag4.14
$$
The equation \thetag{4.14} is integrable in quadratures:
$$
\int\limits^{\,\,\theta}_{\!\!0}\frac{C\,|\sin\theta|\,d\theta}{A(C\,
|\sin\theta|)\,\sin\theta}=t+\tilde C.\hskip -3em
\tag4.15
$$
Hence the whole system of equations \thetag{4.11} for $v(t)$ and
$\theta(t)$ is integrable in quadratures. Its solution contains
two constants of integration $C$ and $\tilde C$.\par
     Further, in order to find trajectories one should integrate
the equation $\dot\bold r=\bold v$. Right hand side of this vectorial
equation in variables $x$, $y$, $v$, $\theta$ is determined by formula
\thetag{4.9} for the velocity vector:
$$
\frac{d}{dt}\Vmatrix x\\ \vspace{1ex} y\endVmatrix=
v\cdot\Vmatrix \cos\theta\\ \vspace{1ex}\sin\theta\endVmatrix.
$$
Knowing $v(t)$ and $\theta(t)$, we find $x(t)$ and $y(t)$ simply by
integration:
$$
\xalignat 2
&\qquad x(t)=x_0+\int\limits^{\,\,t}_{\!\!0}\!\!v(t)\,\cos(\theta(t))\,dt,
&&y(t)=y_0+\int\limits^{\,\,t}_{\!\!0}\!\!v(t)\,\sin(\theta(t))\,dt.
\hskip -3em
\tag4.16
\endxalignat
$$\par
    Function $A(v)$ in formula \thetag{4.7} is arbitrary one. If we
choose it being constant $A(v)=A_0=\const$, then integral in \thetag{4.15}
can be found in explicit form:
$$
\frac{C\,\theta}{A_0}\,\sign(\sin\theta)=t+\tilde C.\hskip -3em
\tag4.17
$$
Let's set up the following Cauchy problem for the equations \thetag{4.11}:
$$
\xalignat 2
&v(t)\,\hbox{\vrule height 8pt depth 8pt width 0.5pt}_{\,t=0}=\,v_0,
&&\theta(t)\,\hbox{\vrule height 8pt depth 8pt width 0.5pt}_{\,t=0}
=\,\theta_0.\hskip -3em
\tag4.18
\endxalignat
$$
Suppose that $A_0>0$ and $\theta_0>0$. Then we have 
$$
\cases
\theta(t)=\theta_0+\omega\,t,\\
\vspace{1ex}
v(t)=\dfrac{A_0}{\omega}\,\sin(\theta_0+\omega\,t).
\endcases\hskip -3em
\tag4.19
$$
Solution of Cauchy problem \thetag{4.18} is determined by formulas
\thetag{4.19} only within some interval of values of parameter $t$:
$$
-\frac{\theta_0}{\omega}\leqslant t\leqslant\frac{\pi-\theta_0}{\omega},
\hskip -3em
\tag4.20
$$
angular velocity $\omega$ in \thetag{4.19} and \thetag{4.20} being
determined by formula
$$
\pagebreak
\omega=\frac{A_0}{v_0}\,\sin\theta_0.\hskip -3em
\tag4.21
$$
If we substitute \thetag{4.19} into formulas \thetag{4.16}, we can
calculate integrals \thetag{4.16} in explicit form. For $x(t)$ and
$y(t)$ this yields
$$
\aligned
&x(t)=x_0-\frac{A_0}{4\,\omega^2}\,\bigl(\cos(2\,\theta_0+2\,\omega\,t)
-\cos(2\,\theta_0)\bigr),\\
\vspace{2ex}
&y(t)=y_0+\frac{A_0\,t}{2\,\omega}-\frac{A_0}{4\,\omega^2}\,\bigl(\sin(2\,
\theta_0+2\,\omega\,t)-\sin(2\,\theta_0)\bigr).
\endaligned\hskip -3em
\tag4.22
$$
From formulas \thetag{4.22} we see that if $A(v)=A_0=\const$,
trajectories of Newtonian dynamical system with force field
\thetag{4.7} are {\bf cycloids}. These are curves drown by a
point of a circle rolling along the straight line parallel to
$OY$ axis.\enddemo
\demo{Example 2} {\bf Spatially non-homogeneous force field with
marked point}. Let's mark the point coinciding with the origin of
Cartesian coordinates and let's direct the vector field $\bold m$
in \thetag{3.1} along the radius-vector:
$$
\xalignat 2
&\quad\bold m=\frac{1}{|\bold r|}\cdot\Vmatrix x\\ \vspace{1ex}y
\endVmatrix,
&&\tilde\bold m=\frac{1}{|\bold r|}\cdot\Vmatrix -y\\ \vspace{1ex}x
\endVmatrix.\hskip -3em
\tag4.23
\endxalignat
$$
Here $|\bold r|=\sqrt{x^2+y^2\,}$. Using explicit formula \thetag{4.23}
for vectors $\bold m$ and $\tilde\bold m$, we can calculate parameters
$\omega_1$ and $\omega_2$ in formulas \thetag{3.3}:
$$
\omega_i=\frac{\tilde m_i}{|\bold r|},\quad i=1,2.\hskip -3em
\tag4.24
$$
Further we use normality equation \thetag{3.18}. When we take into
account \thetag{4.23} and \thetag{4.24}, this equation is written as
$$
\gather
\frac{(x\,A_y-y\,A_x-x\,A_{\theta x}-y\,A_{\theta y}-A_\theta)
\,\cos\theta}{|\bold r|}\,-\hskip -3em\\
\vspace{1.5ex}
-\,\frac{(x\,A_x+y\,A_y-y\,A_{\theta x}+x\,A_{\theta y}
-A_{\theta\theta})\,\sin\theta}{|\bold r|}\,-\hskip -3em
\tag4.25\\
\vspace{1.5ex}
-\,\frac{A\,A_{\theta v}}{v}+\frac{A\,A_\theta}{v^2}
+\frac{A_\theta\,A_{\theta\theta}}{v^2}+\frac{A_\theta\,A_v}{v}=0.
\hskip -3em
\endgather
$$
With the aim of further simplification of the equation \thetag{4.25}
we use polar coordinates in configuration space $M=\Bbb R^2$. Let's
do the following change of variables:
$$
\cases x=\rho\,\cos\gamma,\\ \vspace{0.2ex} y=\rho\,\sin\gamma.\endcases
\hskip -3em
\tag4.26
$$
Polar coordinates $v$ and $\theta$ in the space of velocities are not
changed thereby. From \thetag{4.26} we derive formulas for derivatives:
$$
\aligned
&\frac{\partial}{\partial x}=\cos\gamma\,\frac{\partial}{\partial\rho}-
\frac{\sin\gamma}{\rho}\,\frac{\partial}{\partial\gamma},\\
\vspace{1ex}
&\frac{\partial}{\partial x}=\sin\gamma\,\frac{\partial}{\partial\rho}+
\frac{\cos\gamma}{\rho}\,\frac{\partial}{\partial\gamma}.
\endaligned\hskip -3em
\tag4.27
$$
Taking into account the relationships \thetag{4.26} and \thetag{4.27},
we can write \thetag{4.25} as
$$
\gathered
\frac{(A_\gamma-A_\theta-\rho\,A_{\theta\rho})\,\cos\theta}{\rho}
+\frac{A_\theta\,A_{\theta\theta}}{v^2}+\frac{A_\theta\,A_v}{v}=\\
\vspace{1.5ex}
=\frac{(\rho\,A_\rho+A_{\theta\gamma}-A_{\theta\theta})\,\sin\theta}
{\rho}+\frac{A\,A_{\theta v}}{v}-\frac{A\,A_\theta}{v^2}.
\endgathered\hskip -3em
\tag4.28
$$
If we restrict ourselves with functions $A=A(\rho,v,\theta)$, which do
not depend on $\gamma$, then we easily find the following solution of
the equation \thetag{4.28}:
$$
A=\frac{A(v)\,\cos\theta}{\rho}.\hskip -3em
\tag4.29
$$
Here $A(v)$ is an arbitrary smooth function of one variable. Applying
\thetag{4.2}, we see that the solution \thetag{4.29} corresponds
to Newtonian dynamical system with force field 
$$
\bold F=\frac{A(v)\cdot(\cos\theta\cdot\bold N+\sin\theta\cdot\bold M)}
{\rho}.\hskip -3em
\tag4.30
$$
In Cartesian coordinates in $\Bbb R^2$ force field \thetag{4.30} takes
the following form:
$$
\bold F=A(|\bold v|)\cdot\frac{2\,\bigl<\bold N,\,\bold r\bigr>\cdot\bold N
-\bold r}{|\bold r|^2}.\hskip -3em
\tag4.31
$$\par
\parshape 14 0cm 12.8cm 0cm 12.8cm 0cm 12.8cm
0cm 12.8cm 4.7cm 8.1cm 4.7cm 8.1cm 4.7cm 8.1cm
4.7cm 8.1cm 4.7cm 8.1cm 4.7cm 8.1cm 4.7cm 8.1cm 4.7cm 8.1cm 4.7cm 8.1cm
0cm 12.8cm
    Fig.~4.2 illustrate the geometry of force field \thetag{4.31}. It is
similar to geometry of the field \thetag{4.7}. \vadjust{\vskip 1.6cm
\special{em:graph ris-e42a.gif}\vskip -1.6cm}Vector of force $\bold F$
forms the angle $\theta$ with vector $\bold N$, and it forms the angle
$2\,\theta$ with radius-vector $\bold r$, which, in turn, forms the angle
$\gamma$ with $OX$ axis. Newtonian dynamical system with force field
\thetag{4.30} admits normal shift of curves in $\Bbb R^2$. Similar to
dynamical system considered in example~1 above, this system is integrable
in quadratures. In order to integrate it let's write the equations of
dynamics \thetag{1.2} in variables $\rho$, $\gamma$, $v$, $\theta$. For
components of radius vector $\bold r$ we have the equalities
$$
\aligned
&\bold r=\rho\cdot\Vmatrix \cos\gamma\\ \vspace{1ex}\sin\gamma
\endVmatrix,\\
&\dot\bold r=\dot\rho\cdot\Vmatrix\cos\gamma\\ \vspace{1ex}
\sin\gamma\endVmatrix+\rho\,\dot\gamma\cdot\Vmatrix -\sin\gamma\\
\vspace{1ex}\cos\gamma\endVmatrix.
\endaligned\hskip -2em
\tag4.32
$$
Let's write analogous equality for components of velocity vector, using
Fig~4.2 for this purpose:
$$
\bold v=v\cdot\Vmatrix \cos(\gamma+\theta)\\ \vspace{1ex}\sin(\gamma
+\theta)\endVmatrix=v\,\cos\theta\cdot\Vmatrix\cos\gamma\\ \vspace{1ex}
\sin\gamma\endVmatrix+v\,\sin\theta\cdot\Vmatrix -\sin\gamma\\
\vspace{1ex}\cos\gamma\endVmatrix.
$$
Comparing this formula with \thetag{4.32}, from $\dot\bold r=\bold v$ we
derive two equations:
$$
\cases
\dot\rho=v\,\cos\theta,\\ \vspace{1ex}
\rho\,\dot\gamma=v\,\sin\theta.
\endcases\hskip -3em
\tag4.33
$$
For the derivative of velocity vector we have the expansion, which
follows from the fact that vector $\bold v$ forms the angle $\gamma
+\theta$ with $OX$ axis:
$$
\dot\bold v=\dot v\cdot\Vmatrix \cos(\gamma+\theta)\\ \vspace{1ex}
\sin(\gamma+\theta)\endVmatrix+v\,(\dot\gamma+\dot\theta)
\cdot\Vmatrix -\sin(\gamma+\theta)\\ \vspace{1ex}\cos(\gamma+\theta)
\endVmatrix.
$$
Force vector $\bold F$ in \thetag{4.30} forms the angle $2\,\theta+\gamma
=\theta+(\gamma+\theta)$ with $OX$ axis (see Fig\.~~4.2), its modulus
is equal  to  $A(v)/\rho$.  Using  elementary  identities  from 
trigonometry,
for the vector $\bold F$ we obtain the expansion 
$$
\bold F=\frac{A(v)}{\rho}\,\cos\theta\cdot\Vmatrix \cos(\gamma+\theta)\\
\vspace{1ex}\sin(\gamma+\theta)\endVmatrix+\frac{A(v)}{\rho}\,\sin\theta
\cdot\Vmatrix -\sin(\gamma+\theta)\\ \vspace{1ex}\cos(\gamma+\theta)
\endVmatrix.
$$
In essential, it coincides with \thetag{4.30}. Now from
$\dot\bold v=\bold F$ we derive two more differential equations in
addition to the equation \thetag{4.33}:
$$
\cases
\dot v=\dfrac{A(v)}{\rho}\,\cos\theta,\\ \vspace{1ex}
v\,(\dot\gamma+\dot\theta)=\dfrac{A(v)}{\rho}\,\sin\theta.
\endcases\hskip -3em
\tag4.34
$$
Let's express $\dot\gamma$ through $\rho$, $v$, and $\theta$ by means
of \thetag{4.33} and substitute the obtained expression into \thetag{4.34}.
As a result the equations \thetag{4.34} turn to 
$$
\cases
\dot v=\dfrac{A(v)}{\rho}\,\cos\theta,\\ \vspace{1ex}
v\,\dot\theta=\left(\!\dfrac{A(v)}{\rho}-\dfrac{v^2}{\rho}\right)
\,\sin\theta.\endcases\hskip -3em
\tag4.35
$$
Similar to \thetag{4.11}, the equations \thetag{4.35} can be treated
as differential equations determining function $v=v(\theta)$ in
parametric form:
$$
\frac{dv}{d\theta}=\frac{\dot v}{\vphantom{V^{V^V}}\dot\theta}=
\frac{v\,A(v)}{A(v)-v^2}\,\cotan\theta.\hskip -3em
\tag4.36
$$
Differential equation \thetag{4.36} can be integrated in quadratures:
$$
\int\limits^{\,\,v}_{\!\!0}\frac{A(v)-v^2}{v\,A(v)}\,dv=\ln(|\sin\theta|)
+\const.\hskip -3em
\tag4.37
$$
Having determined the function $v(\theta)$ from \thetag{4.37}, from
\thetag{4.33} we derive the differential equation for the function
$\rho(\theta)$:
$$
\frac{d\rho}{d\theta}=\frac{\dot\rho}{\vphantom{V^{V^V}}\dot\theta}
=\frac{\rho\,v^2}{A(v)-v^2}\,\cotan\theta.\hskip -3em
\tag4.38
$$
The equation \thetag{4.38}, which we have just obtained, is also
integrable in quadratures:
$$
\ln\rho=\int\limits^{\,\,\theta}_{\!\!0}\frac{v(\theta)^2\,\cotan\theta}
{A(v(\theta))-v(\theta)^2}\,d\theta+\const.\hskip -3em
\tag4.39
$$
Let's determine the function $\rho(\theta)$ from \thetag{4.39} and then
substitute it into second equation \thetag{4.35} together with the
function $v(\theta)$ obtained from \thetag{4.37}:
$$
\dot\theta=\frac{A(v)-v^2}{\rho\,v}\,\sin\theta.\hskip -3em
\tag4.40
$$
The solution of \thetag{4.40} is the function $\theta=\theta(t)$,
which is calculated in quadratures according th the following
formula:
$$
\int\limits^{\,\,\theta}_{\!\!0}\frac{\rho(\theta)\,v(\theta)}
{A(v(\theta))-v(\theta)^2}\,\frac{d\theta}{\sin\theta}=t+\const.
\hskip -3em
\tag4.41
$$
Having determined $\theta=\theta(t)$ from \thetag{4.41} and having
substituted it into $v(\theta)$ and $\rho(\theta)$, we get the
functions $v=v(t)$ and $\rho=\rho(t)$. We are to determine the last
of four functions $v(t)$, $\rho(t)$, $\theta(t)$, $\gamma(t)$. It
can be determined from second equation \thetag{4.33} as a result of
simple integration:
$$
\gamma(t)=\int\limits^{\,\,t}_{\!\!0}\frac{v(t)\,\sin(\theta(t))}
{\rho(t)}\,dt+\const.\hskip -3em
\tag4.42
$$
Formulas \thetag{4.37}, \thetag{4.39}, \thetag{4.41}, and \thetag{4.42}
prove that Newtonian dynamical system with force field \thetag{4.31}
is integrable in quadratures.\enddemo
\head
\S\,5. Special classes of solutions of normality equations.
\endhead
     Non-trivial examples of dynamical systems admitting the normal
shift, which are considered in \S\,4, were found in paper \cite{Bol3}
(see also preprint \cite{Bol2}). Their existence stimulated further
investigations and caused the rise of the theory of Newtonian
dynamical systems admitting the normal shift on arbitrary Riemannian
and Finslerian manifolds. Thereby some special classes of such systems
were found.\par
     {\bf 1. Geodesic flows of conformally euclidean metrics}.
Transferring classical Bonnet construction to Riemannian geometry
one obtains the construction of geodesic normal shift, which reproduce
all properties of initial construction. In particular this means that
geodesic flow is a Newtonian dynamical system admitting the normal
shift of hypersurfaces in the metric of that manifold, where it is
determined. The space $M=\Bbb R^2$ turns to Riemannian manifold, if
we define conformally Euclidean metric in it with the following
components:
$$
g_{ij}=e^{-2f}\,\delta_{ij}\text{, \ where \ }\delta_{ij}=\cases 1
\text{\ \ for \ }i=j,\\0\text{\ \ for \ }i\neq j.\endcases\hskip -3em
\tag5.1
$$
Here $f=f(\bold r)$ is some scalar field determining conformal
factor $e^{-2f}$ in metric \thetag{5.1}. Geodesic flow of
conformally Euclidean metric \thetag{5.1} in standard metric of
$\Bbb R^2$ is represented by Newtonian dynamical system with
force field
$$
\bold F=-|\bold v|^2\cdot\nabla f+2\,\bigl<\nabla f,\,\bold v\bigr>
\cdot\bold v.\hskip -2em
\tag5.2
$$
Conformal factor $e^{-2f}$  changes  lengths  of  vectors  when 
transferring
from standard  Euclidean metric in $\Bbb R^2$ to the metric \thetag{5.1}.
However, it doesn't change the measure of angles between vectors. For
this reason force fields \thetag{5.2} admit the normal shift of curves
in standard Euclidean metric in $\Bbb R^2$.\par
     {\bf 2. Metrizable dynamical systems}. Problem of metrizability
arises in connection with force fields of the form \thetag{5.2}. The
matter is that the set of trajectories of one Newtonian dynamical
system with force field $\bold F_1$ can be immersed into the set of
trajectories of other dynamical system with the force field $\bold F_2$.
In such situation we say that first system inherits trajectories of
second one. Dynamical systems that inherit trajectories of the system
with force field \thetag{5.2} were called {\bf metrizable systems}.
The problem of metrizability consists in describing all {\it metrizable
Newtonian dynamical systems admitting the normal shift}. This problem was
investigated and was completely solved in papers \cite{Shr1} and \cite{Shr3}.
Force fields of metrizable dynamical systems admitting the normal shift
of curves in $\Bbb R^2$ have the following form
$$
\bold F(\bold r,\bold v)=-|\bold v|^2\cdot\nabla f+2\,\bigl<\nabla f,\,
\bold v\bigr>\cdot\bold v+\frac{\bold v}{|\bold v|}\cdot
H(|\bold v|\,e^{-f})\,e^f\!.\hskip -1em
\tag5.3
$$
Force fields \thetag{5.3} are parameterized by two smooth scalar
functions: the function of two variables $f=f(x,y)$ and the function
of one variable $H=H(v)$.\par
     {\bf 3. Dynamical systems of multidimensional type}. Force fields
of the form \thetag{5.3} are trivial in some sense: normal shift of
curves along trajectories of Newtonian dynamical systems with such
force fields coincides with geodesic normal shift defined by metric of
the form \thetag{5.1}. Therefore the problem of constructing
non-metrizable dynamical systems admitting the normal shift was
actual. In two-dimensional $n=2$ case this problem is solved by the
examples~1 and 2 considered in \S\,4. Some examples of non-metrizable
dynamical systems admitting the normal shift in multidimensional case
$n\geqslant 3$ were constructed in paper \cite{Bro1}. Then we had the
problem of exhaustive description of all solutions of complete
system of normality equations in multidimensional case. This problem
was considered in paper \cite{Bol9}. However, \pagebreak there an error
in calculations were committed. This error was corrected in \cite{Shr5}.
As a result the following theorem was proved.
\proclaim{Theorem 5.1} Newtonian dynamical system in Riemannian
manifold $M$ of the dimension $n\geqslant 3$ admits normal shift
if and if its force field $\bold F$ is given by formula
$$
\bold F=\frac{h(W)\cdot\bold N}{W_v}-|\bold v|\cdot\frac{2\,
\bigl<\nabla W,\,\bold N\bigr>\cdot\bold N-\nabla W}{W_v},
\hskip -2em
\tag5.4
$$
where $W$ is arbitrary function on $M$ depending on additional
numeric parameter $v=|\bold v|$ and satisfying the condition
$W_v=\partial W/\partial v\neq 0$, while $h=h(v)$ is an arbitrary
function of one variable.
\endproclaim
    In two-dimensional case $n=2$ force fields of the form
\thetag{5.4} also satisfy the normality equations (see
theorem~5.2 below). Corresponding dynamical systems in $M=\Bbb R^2$
are called systems {\bf of multidimensional type}.
\proclaim{Theorem 5.2} Force fields of the form \thetag{5.4} satisfy
the normality equations \thetag{1.16} and \thetag{1.17}.
\endproclaim
\demo{Proof} The equations \thetag{1.16} and \thetag{1.17} are
equivalent to the equations \thetag{1.1}, which they were derived
from. Therefore we can consider the expansion \thetag{1.3} for
$\bold F$ and check the equations \thetag{1.1} for $A$ and $B$.
Coefficients $A$ and $B$ in the expansion \thetag{1.3} are determined
by formulas \thetag{1.5} and \thetag{1.7}.  For the field \thetag{5.4}
they yield
$$
\xalignat 2
&\quad A=\frac{h(W)-\bigl<\nabla W,\,\bold v\bigr>}{W_v},
&&B=|\bold v|\,\frac{\bigl<\nabla W,\,\bold M\bigr>}{W_v}.
\hskip -2em
\tag5.5
\endxalignat
$$
The function $W=W(\bold r,v)$ depends on velocity vector $\bold v$
only through its dependence on $v=|\bold v|$. Therefore $\tilde\nabla W=
W_v\cdot\bold N$. Using this fact, we can explicitly calculate the
velocity gradient for the function $A$:
$$
\tilde\nabla A=\left(\!h'(W)-\frac{\bigl<\nabla W_v,\,\bold v\bigr>}{W_v}
-\frac{h(W)-\bigl<\nabla W,\,\bold v\bigr>}{(W_v)^2}\,W_{vv}\!\right)
\cdot\bold N-\frac{\nabla W}{W_v}.
$$
Then we find the coefficients $\alpha_3$ and $\alpha_4$ in the expansions
\thetag{1.4}:
$$
\gather
\alpha_4=\frac{\bigl<\nabla W,\,\bold M\bigr>}{W_v}.\hskip -2em
\tag5.6\\
\vspace{1ex}
\quad\alpha_3=\frac{\bigl<\nabla W,\,\bold v\bigr>}{(W_v)^2}\,W_{vv}
-\frac{h(W)}{(W_v)^2}\,W_{vv}
+h'(W)-\frac{\bigl<\nabla W_v,\,\bold v\bigr>}{W_v}
-\frac{\bigl<\nabla W,\,\bold N\bigr>}{W_v}.
\hskip -2em
\tag5.7
\endgather
$$
Let's compare \thetag{5.6} with formula for $B$ in \thetag{5.5}
and get sure that first normality equation \thetag{1.1}
for the force field \thetag{5.4} \pagebreak is fulfilled.\par
    Now let's check second normality equation in \thetag{1.1}.
First let's calculate velocity gradient of the function $B$:
$$
\tilde\nabla B=\frac{\bigl<\nabla W,\,\bold M\bigr>}{W_v}\cdot\bold N
+\frac{\bigl<\nabla W_v,\,\bold M\bigr>}{W_v}\cdot\bold v
-\frac{\bigl<\nabla W,\,\bold M\bigr>}{(W_v)^2}\,W_{vv}\cdot\bold v
-\frac{\bigl<\nabla W,\,\bold N\bigr>}{W_v}\cdot\bold M.
$$
Then we find coefficients $\beta_3$ and $\beta_4$ in the expansions
\thetag{1.4}:
$$
\gather
\beta_4=-\frac{\bigl<\nabla W,\,\bold N\bigr>}{W_v}.\hskip -2em
\tag5.8\\
\vspace{1ex}
\beta_3=\frac{\bigl<\nabla W,\,\bold M\bigr>}{W_v}
+|\bold v|\,\frac{\bigl<\nabla W_v,\,\bold M\bigr>}{W_v}
-|\bold v|\,\frac{\bigl<\nabla W,\,\bold M\bigr>}{(W_v)^2}\,W_{vv}.
\hskip -2em
\tag5.9
\endgather
$$
Coefficients $\alpha_2$ and $\beta_1$ in the expansions \thetag{1.4}
are determined as follows:
$$
\xalignat 2
&\alpha_2=\bigl<\nabla A,\,\bold M\bigr>,
&&\beta_1=\bigl<\nabla B,\,\bold N\bigr>.
\endxalignat
$$
By direct calculations according to these formulas we get:
$$
\align
&\gathered
\alpha_2=\frac{h'(W)}{W_v}\,\bigl<\nabla W,\,\bold M\bigr>
-\frac{h(W)}{(W_v)^2}\,\bigl<\nabla W_v,\,\bold M\bigr>\,-\\
\vspace{1ex}
-\,\frac{\nabla_{\bold M}\nabla_{\bold v}W}{W_v}
+\frac{\bigl<\nabla W,\,\bold v\bigr>\,\bigl<\nabla W_v,\,\bold M\bigr>}
{(W_v)^2},
\endgathered\hskip -3em
\tag5.10\\
\vspace{2ex}
&\beta_1=\frac{\nabla_{\bold v}\nabla_{\bold M}W}{W_v}
-\frac{\bigl<\nabla W,\,\bold M\bigr>\,\bigl<\nabla W_v,\,\bold v\bigr>}
{(W_v)^2}.\hskip -3em
\tag5.11
\endalign
$$
Here for the sake of brevity we have introduced the following notations:
$$
\nabla_{\bold M}\nabla_{\bold v}W=\nabla_{\bold v}\nabla_{\bold M}W=
\sum^2_{i=1}\sum^2_{j=1}M^i\,v^j\,\nabla_i\nabla_jW.
$$
Adding the equalities \thetag{5.10} and \thetag{5.11} for the
sum  $\alpha_2+\beta_1$ we obtain formula
$$
\gathered
\alpha_2+\beta_1=\frac{h'(W)}{W_v}\,\bigl<\nabla W,\,\bold M\bigr>
-\frac{h(W)}{(W_v)^2}\,\bigl<\nabla W_v,\,\bold M\bigr>\,+\\
\vspace{1ex}
+\,\frac{\bigl<\nabla W,\,\bold v\bigr>\,\bigl<\nabla W_v,\,\bold M\bigr>}
{(W_v)^2}-\frac{\bigl<\nabla W,\,\bold M\bigr>\,\bigl<\nabla W_v,\,
\bold v\bigr>}{(W_v)^2}.
\endgathered\hskip -3em
\tag5.12
$$
Now let's calculate $\alpha_3-\beta_4$. Taking into account \thetag{5.7}
and \thetag{5.8}, we find
$$
\alpha_3-\beta_4=\frac{\bigl<\nabla W,\,\bold v\bigr>-h(W)}{(W_v)^2}
\,W_{vv}+h'(W)-\frac{\bigl<\nabla W_v,\,\bold v\bigr>}{W_v}.\hskip -2em
\tag5.13
$$
Using formula \thetag{5.9}, we calculate another expression
$$
\frac{B}{|\bold v|}-\beta_3=|\bold v|\,\frac{\bigl<\nabla W,\,
\bold M\bigr>}{(W_v)^2}\,W_{vv}-|\bold v|\,\frac{\bigl<\nabla W_v,
\,\bold M\bigr>}{W_v}.\hskip -3em
\tag5.14
$$
Let's multiply the expression \thetag{5.14} by $|\bold v|^{-1}\,A$,
then multiply \thetag{5.13} by $|\bold v|^{-1}\,B$, and afterwards,
let's add the obtained expressions and subtract \thetag{5.12} from
the sum. Thereby let's use formulas \thetag{5.5}. Left hand side of 
resulting equality coincides with left hand side of second normality
equation in \thetag{1.1}. While right hand side is identically zero.
Theorem~5.2 is proved.\qed\enddemo
    Class of force fields of multidimensional type \thetag{5.4}
comprises all previous examples of force fields of two-dimensional
Newtonian dynamical systems admitting the normal shift. Let's do
the following substitution into the formula \thetag{5.4}:
$$
\xalignat 2
&\quad W(\bold r,v)=v\,e^{-f(\bold r)},&&h(v)=H(v).\hskip -3em
\tag5.15
\endxalignat
$$
As a result of substituting \thetag{5.15} into \thetag{5.15} this
formula takes the form \thetag{5.3}. Further choice $H(v)=0$ reduces
it to the form \thetag{5.2}.\par
     Examples~1 and 2 considered in \S\,4, are also obtained by reduction
from formula \thetag{5.4}. Let's determine $a(v)$ as the solution of
ordinary differential equation
$$
a'(v)=\frac{a(v)}{A(v)}\text{, \ where \ }A(v)\neq 0.\hskip -3em
\tag5.16
$$
Then we define the function $f(\bold r)=\bigl<\bold m,\,\bold r\bigr>$,
where $\bold m=\const$, and we do the following substitution into
the formula \thetag{5.4} for $\bold F$:
$$
\xalignat 2
&\quad W(\bold r,v)=a(v)\,e^{-f(\bold r)},&&h(v)=0.\hskip -3em
\tag5.17
\endxalignat
$$
The substitution \thetag{5.17} transforms \thetag{5.4} to the form
\thetag{4.8}. This corresponds to the example~1 in \S\,4 above.\par
     In order to obtain the force field \thetag{4.13} from the
example~2 in \S\,4 we should choose another function $f=f(\bold r)$
in \thetag{5.17}. Let's take $f(\bold r)=|\bold r|$. The function
$a(v)$ will be the solution of differential equation \thetag{5.16}
as before.\par
    Functions of the form \thetag{5.4} exhaust force fields of dynamical
systems admitting the normal shift in arbitrary Riemannian manifolds
of the dimension $n\geqslant 3$. Is it true in two-dimensional case
$n=2$\,? This should be clear from further analysis of normality equations
below.
\head
\S\,6. Point symmetries of reduced normality equations.
\endhead
    Symmetry analysis is one of the powerful tools for constructing
special solutions of differential equations. Let's apply it in order
to construct new classes of Newtonian dynamical systems admitting the
normal shift in $\Bbb R^2$. By means of scalar ansatz \thetag{2.4}
we reduced normality equations \thetag{1.16} and \thetag{1.17}
to one partial differential equation. This is reduced normality
equation, which can be written as \thetag{3.18}, \thetag{4.5}, or
\thetag{4.25}. All these forms of reduced normality equation are
equivalent to each other. For the further analysis we choose the
equation \thetag{4.5} as the most simple:
$$
\align
(A_y&-A_{\theta x})\,\cos\theta-(A_x+A_{\theta y})\,\sin\theta\,+
\hskip -2em\\
\vspace{-1.0ex}
&\tag6.1\\
\vspace{-1.5ex}
&+\,\frac{A\,A_\theta}{v^2}+\frac{A_\theta\,A_{\theta\theta}}{v^2}
+\frac{A_\theta\,A_v}{v}=\frac{A\,A_{\theta v}}{v}.
\hskip -2em
\endalign
$$
Let's consider the space $\Bbb R^5=\Bbb R^4\oplus\Bbb R$, coordinates
of points in which are the independent variables $x$, $y$, $v$, $\theta$
of the equation \thetag{6.1} and the variable $A$. Here $A$ is also
treated as independent variable:
$$
\Vmatrix x\\ y\\ v\\ \theta\\ A\endVmatrix\in TM\oplus\Bbb R=\Bbb R^5.
\hskip -2em
\tag6.2
$$
Vector field $U$ in five-dimensional space $\Bbb R^5$ from \thetag{6.2}
can be defined as linear differential operator of the first order
$$
U=X\,\frac{\partial}{\partial x}+Y\,\frac{\partial}{\partial y}
+V\,\frac{\partial}{\partial v}+\Theta\,\frac{\partial}{\partial
\theta}+a\,\frac{\partial}{\partial A}.\hskip -2em
\tag6.3
$$
Coefficients $X$, $Y$, $V$, $\Theta$, and $a$ in differential operator
\thetag{6.3} are the functions of $x$, $y$, $v$, $\theta$, $A$.
Vector field \thetag{6.3} corresponds to the local one-parametric group
of local diffeomorphisms in $\Bbb R^5$ (see \cite{Kob1}):
$$
\varphi_t\!:\Vmatrix x\\ \vspace{0.5ex}y\\  \vspace{0.5ex}v\\
\vspace{0.5ex}\theta\\ \vspace{0.5ex}A\endVmatrix\mapsto
\Vmatrix\tilde x\\ \vspace{0.5ex}\tilde y\\ \vspace{0.5ex}\tilde v\\
\vspace{0.5ex}\tilde\theta\\ \vspace{0.5ex}\tilde A\endVmatrix=
\Vmatrix\tilde x(x,y,v,\theta,A,t)\\ \vspace{0.5ex}\tilde y(x,y,v,\theta,
A,t)\\ \vspace{0.5ex}\tilde v(x,y,v,\theta,A,t)\\ \vspace{0.5ex}\tilde
\theta(x,y,v,\theta,A,t)\\ \vspace{0.5ex}\tilde A(x,y,v,\theta,A,t)
\endVmatrix.\hskip -2em
\tag6.4
$$
Functions $\tilde x$, $\tilde y$, $\tilde v$, $\tilde\theta$, $\tilde A$
determining the transformation \thetag{6.4} should satisfy the system of
differential equations with respect to parameter $t$:
$$
\pagebreak
\cases
\dot{\tilde x}=X(\tilde x,\tilde y,\tilde v,\tilde\theta,\tilde A),\\
\dot{\tilde y}=Y(\tilde x,\tilde y,\tilde v,\tilde\theta,\tilde A),\\
\dot{\tilde v}=V(\tilde x,\tilde y,\tilde v,\tilde\theta,\tilde A),\\
\dot{\tilde\theta}=\Theta(\tilde x,\tilde y,\tilde v,\tilde\theta,
\tilde A),\\
\dot{\tilde A}=a(\tilde x,\tilde y,\tilde v,\tilde\theta,\tilde A).
\endcases\hskip -2em
\tag6.5
$$
For $t=0$ the map $\varphi_t$ in \thetag{6.4} is the identical map.
This expressed by the conditions
$$
\xalignat 5
&\qquad\tilde x\,\hbox{\vrule height 8pt depth 8pt width 0.5pt}_{\,t=0}=x,
&&\tilde y\,\hbox{\vrule height 8pt depth 8pt width 0.5pt}_{\,t=0}=y,
&&\tilde v\,\hbox{\vrule height 8pt depth 8pt width 0.5pt}_{\,t=0}=v,
&&\tilde\theta\,\hbox{\vrule height 8pt depth 8pt width 0.5pt}_{\,t=0}
=\,\theta,
&&\tilde A\,\hbox{\vrule height 8pt depth 8pt width 0.5pt}_{\,t=0}=\,A.
\hskip -1em
\tag6.6
\endxalignat
$$
Assuming vector field $U$ in \thetag{6.3} to be smooth, we consider
Taylor expansions of the functions $\tilde x$, $\tilde y$, $\tilde v$,
$\tilde\theta$, $\tilde A$ from \thetag{6.4} with respect to parameter
$t$ in the neighborhood of the point $t=0$. Thereby we restrict ourselves
by the terms of the first order. From the relationships \thetag{6.5}
and \thetag{6.6} we derive:
$$
\align
&\tilde x=x+X(x,y,v,\theta,A)\cdot t+O(t^2),\hskip -2em\\
&\tilde y=y+Y(x,y,v,\theta,A)\cdot t+O(t^2),\hskip -2em\\
&\tilde v=v+V(x,y,v,\theta,A)\cdot t+O(t^2),\hskip -2em
\tag6.7\\
&\tilde\theta=\theta+\Theta(x,y,v,\theta,A)\cdot t+O(t^2),
\hskip -2em\\
&\tilde A=A+a(x,y,v,\theta,A)\cdot t+O(t^2).\hskip -2em
\endalign
$$
In order to invert the transformations forming one-parametric group
in $t$ we should only change $t$ by $-t$, i\,~e\. $(\varphi_t)^{-1}
=\varphi_{-t}$. Therefore transformations inverse to \thetag{6.7}
can be written as follows:
$$
\varphi_{-t}\!:\Vmatrix x\\ \vspace{0.5ex}y\\  \vspace{0.5ex}v\\
\vspace{0.5ex}\theta\\ \vspace{0.5ex}A\endVmatrix\mapsto
\Vmatrix\tilde x\\ \vspace{0.5ex}\tilde y\\ \vspace{0.5ex}\tilde v\\
\vspace{0.5ex}\tilde\theta\\ \vspace{0.5ex}\tilde A\endVmatrix=
\Vmatrix\tilde x(x,y,v,\theta,A,-t)\\ \vspace{0.5ex}\tilde y(x,y,v,\theta,
A,-t)\\ \vspace{0.5ex}\tilde v(x,y,v,\theta,A,-t)\\ \vspace{0.5ex}\tilde
\theta(x,y,v,\theta,A,-t)\\ \vspace{0.5ex}\tilde A(x,y,v,\theta,A,-t)
\endVmatrix.\hskip -2em
\tag6.8
$$
The expansions analogous to \thetag{6.7} for transformations \thetag{6.8}
are written as
$$
\align
&\tilde x=x-X(x,y,v,\theta,A)\cdot t+O(t^2),\hskip -2em\\
&\tilde y=y-Y(x,y,v,\theta,A)\cdot t+O(t^2),\hskip -2em\\
&\tilde v=v-V(x,y,v,\theta,A)\cdot t+O(t^2),\hskip -2em
\tag6.9\\
&\tilde\theta=\theta-\Theta(x,y,v,\theta,A)\cdot t+O(t^2),
\hskip -2em\\
&\tilde A=A-a(x,y,v,\theta,A)\cdot t+O(t^2).\hskip -2em
\endalign
$$\par
    Each function $A=A(x,y,v,\theta)$ determines some hypersurface
$\Gamma$ in five-dimensio\-nal space \thetag{6.2}, $\Gamma$ is a graph
of the function $A=A(x,y,v,\theta)$. The solution of the equation
\thetag{6.1} is exactly that function we need. Applying transformation
\thetag{6.7}, we get new hypersurface $\Gamma_t=\varphi_t(\Gamma)$.
For sufficiently small value of parameter $t\to 0$ hypersurface
$\Gamma_t$ slightly differs from initial hypersurface $\Gamma$. The
hypersurface $\Gamma_t$ or some part of this hypersurface can be
considered as a graph for some other function $A=A^t(x,y,v,\theta)$.
So the transformations \thetag{6.4} defined by vector field
\thetag{6.3} can be treated as transformations in the set of functions
$A=A(x,y,v,\theta)$.
\definition{Definition 6.1} Vector field \thetag{6.3} is called {\bf the
field of point symmetry} for the equation \thetag{6.1} if corresponding
transformations $\varphi_t$ in \thetag{6.4} transforms each solution of
this equation into another solution of the same equation.
\enddefinition
\noindent Field of point symmetry of the equation \thetag{6.1} form
Lie subalgebra in the Lie algebra of all vector fields in five-dimensional
space \thetag{6.2}. Denote it by $\goth L$. Algebra $\goth L$ is called
the algebra of point symmetries of the equation \thetag{6.1}. Presently
there is a well-developed theory for calculating point symmetry algebra
of arbitrary differential equation (see more details in \cite{Ibr1} and
\cite{Olv1}). Let's apply it to the equation \thetag{6.1}.\par
     Let $A(x,y,v,\theta)$ be the solution of the equation \thetag{6.1}.
We substitute the function $A=A(x,y,v,\theta)$ into the arguments of the
functions $\tilde x$, $\tilde y$, $\tilde v$, $\tilde\theta$, $\tilde A$,
which determine the transformation \thetag{6.4}. As a result we get five
functions
$$
\align
&\tilde x=\tilde x(x,y,v,\theta,A(x,y,v,\theta),t),\hskip -3em\\
&\tilde y=\tilde y(x,y,v,\theta,A(x,y,v,\theta),t),\hskip -3em\\
&\tilde v=\tilde v(x,y,v,\theta,A(x,y,v,\theta),t),\hskip -3em
\tag6.10\\
&\tilde\theta=\tilde\theta(x,y,v,\theta,A(x,y,v,\theta),t),\hskip -3em\\
&\tilde A=\tilde A(x,y,v,\theta,A(x,y,v,\theta),t)\hskip -3em\\
\endalign
$$
of four variables (the variable $t$ here is considered as parameter).
Functions \thetag{6.10} in implicit (parametric) form determine the function
$\tilde A=A^t(\tilde  x,\tilde  y,\tilde  v,\tilde\theta)$.  In 
order to
calculate this function explicitly we should use first four relationships
\thetag{6.10} for to express $x$, $y$, $v$, and $\theta$ through
$\tilde x$, $\tilde y$, $\tilde v$, and $\tilde\theta$. Then the obtained
expressions for $x$, $y$, $v$, and $\theta$ should be substituted into
fifth relationship \thetag{6.10}. Denote by $J$ the Jacoby matrix
of the change of variables given by first four relationships \thetag{6.10}:
$$
J=\frac{\partial(\tilde x,\tilde y,\tilde v,\tilde\theta)}{\partial(x,
y,v,\theta)}.\hskip -3em
\tag6.11
$$
Let's consider some components in Jacoby matrix \thetag{6.11}:
$$
\aligned
&J^1_1=\tilde x_x(x,y,v,\theta,A,t)+\tilde x_A(x,y,v,\theta,A,t)\cdot A_x,\\
&J^1_2=\tilde x_y(x,y,v,\theta,A,t)+\tilde x_A(x,y,v,\theta,A,t)\cdot A_y,\\
&.\ .\ .\ .\ .\ .\ .\ .\ .\ .\ .\ .\ .\ .\ .\ .\ .\ .\ .\ .\ .\ .\ .\ .\ 
.\ .\ .\ .\ .\ .\ .\ .\ .\ .\\
&J^4_3=\tilde\theta_v(x,y,v,\theta,A,t)+\tilde\theta_A(x,y,v,\theta,A,t)
\cdot A_v,\\
&J^4_4=\tilde\theta_\theta(x,y,v,\theta,A,t)+\tilde\theta_A(x,y,v,\theta,
A,t)\cdot A_\theta.
\endaligned\hskip -3em
\tag6.12
$$
From \thetag{6.12} we see that components of Jacoby matrix \thetag{6.11}
depend on variables $x$, $y$, $v$, $\theta$, and $A$. Moreover, they
depend linearly on derivatives $A_x$, $A_y$, $A_v$, and $A_\theta$:
$$
\pagebreak
J^k_q=J^k_q(x,y,v,\theta,A,A_x,A_y,A_v,A_\theta,t).\hskip -3em
\tag6.13
$$
Denote by $I=J^{-1}$ the matrix inverse to $J$. Its components depend
on the same quantities as components of $J$ in \thetag{6.13}:
$$
I^k_q=I^k_q(x,y,v,\theta,A,A_x,A_y,A_v,A_\theta,t).\hskip -3em
\tag6.14
$$
Though the dependence on derivatives $A_x$, $A_y$, $A_v$, and $A_\theta$
in \thetag{6.14} is not linear.\par
    We use components of matrix $I$ for to calculate partial derivatives
of the function $A^t(\tilde x,\tilde y,\tilde v,\tilde\theta)$, which is
defined implicitly by the relationships \thetag{6.10}. Remember that the
variables $x$, $y$, $v$, $\theta$ are expressed through $\tilde x$,
$\tilde y$, $\tilde v$, $\tilde\theta$ due to first four relationships
\thetag{6.10}, and $I$ is Jacoby matrix for such expressions:
$$
I=\frac{\partial(x,y,v,\theta)}{\partial(\tilde x,\tilde y,\tilde v,
\vphantom{V^{V^V}}\tilde\theta)}.\hskip -3em
\tag6.15
$$
Then, due to \thetag{6.15}, for the derivative $A^t_{\tilde x}$ from
\thetag{6.10} we obtain the expression
$$
\aligned
A^t_{\tilde x}&=\bigl(\tilde A_x(x,y,v,\theta,A,t)+\tilde A_A(x,y,v,
\theta,A,t)\cdot A_x\bigr)\cdot I^1_1\,+\\
&+\,\bigl(\tilde A_y(x,y,v,\theta,A,t)+\tilde A_A(x,y,v,
\theta,A,t)\cdot A_y\bigr)\cdot I^2_1\,+\\
&+\,\bigl(\tilde A_v(x,y,v,\theta,A,t)+\tilde A_A(x,y,v,
\theta,A,t)\cdot A_v\bigr)\cdot I^3_1\,+\\
&+\,\bigl(\tilde A_\theta(x,y,v,\theta,A,t)+\tilde A_A(x,y,v,
\theta,A,t)\cdot A_\theta\bigr)\cdot I^4_1.
\endaligned
\tag6.16
$$
We have similar expressions for other derivatives:
$$
\gather
\aligned
A^t_{\tilde y}&=\bigl(\tilde A_x(x,y,v,\theta,A,t)+\tilde A_A(x,y,v,
\theta,A,t)\cdot A_x\bigr)\cdot I^1_2\,+\\
&+\,\bigl(\tilde A_y(x,y,v,\theta,A,t)+\tilde A_A(x,y,v,
\theta,A,t)\cdot A_y\bigr)\cdot I^2_2\,+\\
&+\,\bigl(\tilde A_v(x,y,v,\theta,A,t)+\tilde A_A(x,y,v,
\theta,A,t)\cdot A_v\bigr)\cdot I^3_2\,+\\
&+\,\bigl(\tilde A_\theta(x,y,v,\theta,A,t)+\tilde A_A(x,y,v,
\theta,A,t)\cdot A_\theta\bigr)\cdot I^4_2.
\endaligned\hskip -3em
\tag6.17\\
\vspace{2ex}
\aligned
A^t_{\tilde v}&=\bigl(\tilde A_x(x,y,v,\theta,A,t)+\tilde A_A(x,y,v,
\theta,A,t)\cdot A_x\bigr)\cdot I^1_3\,+\\
&+\,\bigl(\tilde A_y(x,y,v,\theta,A,t)+\tilde A_A(x,y,v,
\theta,A,t)\cdot A_y\bigr)\cdot I^2_3\,+\\
&+\,\bigl(\tilde A_v(x,y,v,\theta,A,t)+\tilde A_A(x,y,v,
\theta,A,t)\cdot A_v\bigr)\cdot I^3_3\,+\\
&+\,\bigl(\tilde A_\theta(x,y,v,\theta,A,t)+\tilde A_A(x,y,v,
\theta,A,t)\cdot A_\theta\bigr)\cdot I^4_3.
\endaligned\hskip -3em
\tag6.18\\
\vspace{2ex}
\aligned
A^t_{\tilde\theta}&=\bigl(\tilde A_x(x,y,v,\theta,A,t)+\tilde A_A(x,y,v,
\theta,A,t)\cdot A_x\bigr)\cdot I^1_4\,+\\
&+\,\bigl(\tilde A_y(x,y,v,\theta,A,t)+\tilde A_A(x,y,v,
\theta,A,t)\cdot A_y\bigr)\cdot I^2_4\,+\\
&+\,\bigl(\tilde A_v(x,y,v,\theta,A,t)+\tilde A_A(x,y,v,
\theta,A,t)\cdot A_v\bigr)\cdot I^3_4\,+\\
&+\,\bigl(\tilde A_\theta(x,y,v,\theta,A,t)+\tilde A_A(x,y,v,
\theta,A,t)\cdot A_\theta\bigr)\cdot I^4_4.
\endaligned\hskip -3em
\tag6.19
\endgather
$$
Let's distract ourselves for a while from particular form of right
hand sides of the equalities \thetag{6.16}, \thetag{6.17}, \thetag{6.18},
and \thetag{6.19}. Let's rewrite them in symbolic form keeping only
the indication of functional dependence of those quantities, which are
present in right hand sides of these equalities:
$$
\aligned
&A^t_{\tilde x}=\tilde A_1(x,y,v,\theta,A,A_x,A_y,A_v,A_\theta,t),\\
&A^t_{\tilde y}=\tilde A_2(x,y,v,\theta,A,A_x,A_y,A_v,A_\theta,t),\\
&A^t_{\tilde v}=\tilde A_3(x,y,v,\theta,A,A_x,A_y,A_v,A_\theta,t),\\
&A^t_{\tilde\theta}=\tilde A_4(x,y,v,\theta,A,A_x,A_y,A_v,A_\theta,t).
\endaligned\hskip -3em
\tag6.20
$$
The functions $\tilde A_1$, $\tilde A_2$, $\tilde A_3$, and $\tilde A_4$
\thetag{6.20} do not depend on particular choice of function $A=A(x,y,v,
\theta)$. The function $A$ (its explicit form) is of importance only
when we substitute its derivatives $A_x$, $A_y$, $A_v$, $A_\theta$ into
the arguments of the functions \thetag{6.20}. If we consider $A$, $A_x$,
$A_y$, $A_v$, $A_\theta$ as independent variables and if we complete
$\tilde x$, $\tilde y$, $\tilde v$, and $\tilde A$ in \thetag{6.4} by
functions \thetag{6.20}, then we obtain the extension of transformations
$\varphi_t$ from $\Bbb R^5$ to $\Bbb R^9$. It is called the {\bf first
extension} of local one-parametric group of local transformations
$\varphi_t$.\par
    Differentiating the relationships \thetag{6.20}, we obtain analogous
formulas for second derivatives $A^t_{\tilde x\tilde x}$, $A^t_{\tilde x
\tilde y}$, $A^t_{\tilde x\tilde v}$, $A^t_{\tilde x\tilde\theta}$,
$A^t_{\tilde y\tilde y}$, $A^t_{\tilde y\tilde v}$, $A^t_{\tilde y\tilde
\theta}$, $A^t_{\tilde v\tilde v}$, $A^t_{\tilde v\tilde\theta}$,
$A^t_{\tilde\theta\tilde\theta}$:
$$
\aligned
&A^t_{\tilde x\tilde x}=\tilde A_{11}(x,y,v,\theta,A,A_x,A_y,A_v,A_\theta,
A_{xx},\ldots,A_{\theta\theta},t),\\
&A^t_{\tilde x\tilde y}=\tilde A_{12}(x,y,v,\theta,A,A_x,A_y,A_v,A_\theta,
A_{xx},\ldots,A_{\theta\theta},t),\\
&A^t_{\tilde x\tilde v}=\tilde A_{13}(x,y,v,\theta,A,A_x,A_y,A_v,A_\theta,
A_{xx},\ldots,A_{\theta\theta},t),\\
&A^t_{\tilde x\tilde\theta}=\tilde A_{14}(x,y,v,\theta,A,A_x,A_y,A_v,A_\theta,
A_{xx},\ldots,A_{\theta\theta},t),\\
&.\ .\ .\ .\ .\ .\ .\ .\ .\ .\ .\ .\ .\ .\ .\ .\ .\ .\ .\ .\ .\ .\ .\ .\ 
.\ .\ .\ .\ .\ .\ .\ .\ .\ .\ .\ .\ .\\
&A^t_{\tilde y\tilde\theta}=\tilde A_{24}(x,y,v,\theta,A,A_x,A_y,A_v,A_\theta,
A_{xx},\ldots,A_{\theta\theta},t),\\
&A^t_{\tilde v\tilde v}=\tilde A_{33}(x,y,v,\theta,A,A_x,A_y,A_v,A_\theta,
A_{xx},\ldots,A_{\theta\theta},t),\\
&A^t_{\tilde v\tilde\theta}=\tilde A_{34}(x,y,v,\theta,A,A_x,A_y,A_v,
A_\theta,A_{xx},\ldots,A_{\theta\theta},t),\\
&A^t_{\tilde\theta\tilde\theta}=\tilde A_{44}(x,y,v,\theta,A,A_x,A_y,A_v,
A_\theta,A_{xx},\ldots,A_{\theta\theta},t).
\endaligned\hskip -3em
\tag6.21
$$
Ten functions \thetag{6.21} define {\bf second extension} of local
one-parametric group of transformations $\varphi_t$. Second extension
$\varphi_t$ acts in the space $\Bbb R^{19}$.\par
     First and second extensions of transformations $\varphi_t$ form
local one-parametric groups of local transformations in the spaces
$\Bbb R^9$ and $\Bbb R^{19}$ respectively. They correspond to some vector
fields in $\Bbb R^9$ and $\Bbb R^{19}$, these vector fields are called
{\bf first and second extensions} of vector field $U$ from \thetag{6.3}.
First extension can be represented by the following differential operator:
$$
\gathered
U=X\,\frac{\partial}{\partial x}+Y\,\frac{\partial}{\partial y}
+V\,\frac{\partial}{\partial v}+\Theta\,\frac{\partial}{\partial
\theta}+a\,\frac{\partial}{\partial A}\,+\\
\vspace{1ex}
+\,a_1\,\frac{\partial}{\partial A_x}+\,a_2\,\frac{\partial}
{\partial A_y}+\,a_3\,\frac{\partial}{\partial A_v}+\,a_4\,
\frac{\partial}{\partial A_\theta}.
\endgathered\hskip -3em
\tag6.22
$$
For the second extension we have similar differential operator, but
with more terms:
$$
\gathered
U=X\,\frac{\partial}{\partial x}+Y\,\frac{\partial}{\partial y}
+V\,\frac{\partial}{\partial v}+\Theta\,\frac{\partial}{\partial
\theta}+a\,\frac{\partial}{\partial A}\,+\\
\vspace{1ex}
+\,a_1\,\frac{\partial}{\partial A_x}+\,a_2\,\frac{\partial}
{\partial A_y}+\,a_3\,\frac{\partial}{\partial A_v}+\,a_4\,
\frac{\partial}{\partial A_\theta}\,+\\
\vspace{1ex}
+\,a_{11}\,\frac{\partial}{\partial A_{xx}}+a_{12}\,\frac{\partial}
{\partial A_{xy}}+a_{13}\,\frac{\partial}{\partial A_{xv}}
+a_{14}\,\frac{\partial}{\partial A_{x\theta}}\,+\\
\vspace{1ex}
+\,a_{22}\,\frac{\partial}{\partial A_{yy}}+a_{23}\,\frac{\partial}
{\partial A_{yv}}+a_{24}\,\frac{\partial}{\partial A_{y\theta}}
+a_{33}\,\frac{\partial}{\partial A_{vv}}+a_{34}\,\frac{\partial}
{\partial A_{v\theta}}+a_{44}\,\frac{\partial}{\partial A_{\theta\theta}}.
\endgathered\hskip -1em
\tag6.23
$$
Coefficients $a_1$, $a_2$, $a_3$, $a_4$, $a_{11}$, $a_{12}$, $a_{13}$,
$a_{14}$, $a_{22}$, $a_{23}$, $a_{24}$, $a_{33}$, $a_{34}$, and $a_{44}$
in \thetag{6.22} and in \thetag{6.23} can be calculated in explicit form.
For to calculate $a_1$, $a_2$, $a_3$, $a_4$ one should consider Taylor
expansions of the functions \thetag{6.17}, \thetag{6.18}, \thetag{6.19},
and \thetag{6.20} with respect to parameter $t$ in the neighborhood
of the point $t=0$:
$$
\aligned
&A^t_{\tilde x}=A_x+a_1(x,y,v,\theta,A,A_x,A_y,A_v,A_\theta)
\cdot t+O(t^2),\\
&A^t_{\tilde y}=A_y+a_2(x,y,v,\theta,A,A_x,A_y,A_v,A_\theta)
\cdot t+O(t^2),\\
&A^t_{\tilde v}=A_v+a_3(x,y,v,\theta,A,A_x,A_y,A_v,A_\theta)
\cdot t+O(t^2),\\
&A^t_{\tilde\theta}=A_x+a_4(x,y,v,\theta,A,A_x,A_y,A_v,A_\theta)
\cdot t+O(t^2).
\endaligned\hskip -3em
\tag6.24
$$
Detailed calculations leading to \thetag{6.24} are standard. They
are not of particular interest here. We shall give only the result
of these calculations. For the function $a_1(x,y,v,\theta,A,A_x,A_y,
A_v,A_\theta)$ we get
$$
\gathered
a_1=a_x+a_A\,A_x-A_x\,X_x-A_y\,Y_x-A_v\,V_x-A_\theta\,\Theta_x-\\
\vspace{0.4ex}
-A_x\,X_A\,A_x-A_y\,Y_A\,A_x-A_v\,V_A\,A_x-A_\theta\,\Theta_A\,A_x.
\endgathered\hskip -1em
\tag6.25
$$
Similar expressions are available for $a_2$, $a_3$, and $a_4$:
$$
\align
&\gathered
a_2=a_y+a_A\,A_y-A_x\,X_y-A_y\,Y_y-A_v\,V_y-A_\theta\,\Theta_y-\\
\vspace{0.4ex}
-A_x\,X_A\,A_y-A_y\,Y_A\,A_y-A_v\,V_A\,A_y-A_\theta\,\Theta_A\,A_y.
\endgathered\hskip -3em
\tag6.26\\
\vspace{2ex}
&\gathered
a_3=a_v+a_A\,A_v-A_x\,X_v-A_y\,Y_v-A_v\,V_v-A_\theta\,\Theta_v-\\
\vspace{0.4ex}
-A_x\,X_A\,A_v-A_y\,Y_A\,A_v-A_v\,V_A\,A_v-A_\theta\,\Theta_A\,A_v.
\endgathered\hskip -3em
\tag6.27\\
\vspace{2ex}
&\gathered
a_3=a_\theta+a_A\,A_\theta-A_x\,X_\theta-A_y\,Y_\theta-A_v\,V_\theta
-A_\theta\,\Theta_\theta-\\
\vspace{0.4ex}
-A_x\,X_A\,A_\theta-A_y\,Y_A\,A_\theta-A_v\,V_A\,A_\theta-A_\theta\,
\Theta_A\,A_\theta.
\endgathered\hskip -3em
\tag6.28
\endalign
$$
Formulas for coefficients $a_{11}$, $a_{12}$, $a_{13}$, $a_{14}$,
$a_{22}$, $a_{23}$, $a_{24}$, $a_{33}$, $a_{34}$, and $a_{44}$ in
\thetag{6.23} are more huge. In order to get them one should
consider Taylor expansion of the functions \thetag{6.21} at the
point $t=0$ and should take linear terms in such expansions:
$$
\aligned
&A^t_{\tilde x\tilde x}=A_{xx}+a_{11}(x,y,v,\theta,A,\ldots,A_{v\theta},
A_{\theta\theta})\cdot t+O(t^2),\\
&A^t_{\tilde x\tilde y}=A_{xy}+a_{12}(x,y,v,\theta,A,\ldots,A_{v\theta},
A_{\theta\theta})\cdot t+O(t^2),\\
&.\ .\ .\ .\ .\ .\ .\ .\ .\ .\ .\ .\ .\ .\ .\ .\ .\ .\ .\ .\ .\ .\ .\ .\ 
.\ .\ .\ .\ .\ .\ .\ .\ .\ .\\
&A^t_{\tilde v\tilde\theta}=A_{v\theta}+a_{34}(x,y,v,\theta,A,\ldots,
A_{v\theta},A_{\theta\theta})\cdot t+O(t^2),\\
&A^t_{\tilde\theta\tilde\theta}=A_{\theta\theta}+a_{34}(x,y,v,\theta,A,
\ldots,A_{v\theta},A_{\theta\theta})\cdot t+O(t^2).
\endaligned\hskip -3em
\tag6.29
$$
We do not give explicit formulas for $a_{ij}(x,y,v,\theta,A,\ldots,
A_{v\theta},A_{\theta\theta})$, they are rather huge. We give only
the complete list of arguments in them:
$$
\gathered
x,\,\,y,\,\,v,\,\,\theta,\,A,\,A_x,\,A_y,\,A_v,\,A_\theta,\,A_{xx},
\,A_{xy},\\
\vspace{0.4ex}
A_{xv},\,\,A_{x\theta},\,\,A_{yy},\,\,A_{yv},\,\,A_{y\theta},\,\,
A_{vv},\,\,A_{v\theta},\,\,A_{\theta\theta}.
\endgathered\hskip -3em
\tag6.30
$$\par
     Let $A(x,y,v,\theta)$ be an arbitrary scalar function of
four arguments. Fixing some values of its arguments, we calculate
its value and the values of all its derivatives listed in \thetag{6.30}
for these fixed arguments. As a result we get the point of the space
$\Bbb R^{19}$. Each point of the space $\Bbb R^{19}$ can be obtained
in this way. However, if choosing the function $A(x,y,v,\theta)$ we
restrict our choice to the set of solutions of the normality equation,
we would obtain only the points of some hypersurface in $\Bbb R^{19}$.
The equation \thetag{6.1} itself is the equation of such hypersurface
provided that we treat all quantities in \thetag{6.30} as independent
variables. Remember that reduced normality equation equation \thetag{6.1}
has the following structure:
$$
F(x,y,\ldots,A_{v\theta},A_{\theta\theta})=0.\hskip -3em
\tag6.31
$$
Function $F$ in \thetag{6.31} is determined by left hand side of the
equation \thetag{6.1}. According to general theory of point symmetries
(see \cite{Ibr1} and \cite{Olv1}) vector field \thetag{6.3} is a field of
point symmetry for differential equation \thetag{6.31}, if its second
extension is tangent to hypersurface in $\Bbb R^{19}$ defined by the
equation \thetag{6.31}. This condition is written as 
$$
\gathered
0=X\,\frac{\partial F}{\partial x}+Y\,\frac{\partial F}{\partial y}
+V\,\frac{\partial F}{\partial v}+\Theta\,\frac{\partial F}{\partial
\theta}+a\,\frac{\partial F}{\partial A}
+a_1\,\frac{\partial F}{\partial A_x}+\,a_2\,\frac{\partial F}
{\partial A_y}\,+\\
\vspace{1ex}
+\,a_3\,\frac{\partial F}{\partial A_v}
+a_4\,\frac{\partial F}{\partial A_\theta}
+a_{11}\,\frac{\partial F}{\partial A_{xx}}+a_{12}\,\frac{\partial F}
{\partial A_{xy}}+a_{13}\,\frac{\partial F}{\partial A_{xv}}
+a_{14}\,\frac{\partial F}{\partial A_{x\theta}}\,+\\
\vspace{1ex}
+\,a_{22}\,\frac{\partial F}
{\partial A_{yy}}+a_{23}\,\frac{\partial F}{\partial A_{yv}}
+a_{24}\,\frac{\partial F}{\partial A_{y\theta}}+a_{33}\,\frac{\partial F}
{\partial A_{vv}}+a_{34}\,\frac{\partial F}{\partial A_{v\theta}}
+a_{44}\,\frac{\partial F}{\partial A_{\theta\theta}}.
\endgathered\hskip -1em
\tag6.32
$$
The equality \thetag{6.32} should be fulfilled at all points of the
space $\Bbb R^{19}$, where the equality \thetag{6.31} is fulfilled.
The equality \thetag{6.32} is called the {\bf determining equation}
for the field of point symmetry.
\head
\S\,7. Calculation of the point symmetries.
\endhead
    Let's consider the determining equation \thetag{6.32} for the field
of point symmetry $U$ applied to our case, when the equation \thetag{6.31}
has the form \thetag{6.1}. The equation \thetag{6.1} can be resolved with
respect to one of the second order derivatives:
$$
\align
A_{v\theta}&=\frac{v\,(A_y-A_{\theta x})\,\cos\theta}{A}+\frac{A_\theta
\,A_v}{A}\,-\hskip -2em\\
\vspace{-1.0ex}
&\tag7.1\\
\vspace{-1.5ex}
&-\,\frac{v\,(A_x+A_{\theta y})\,\sin\theta}{A}+\frac{A_\theta}{v}
+\frac{A_\theta\,A_{\theta\theta}}{v\,A}.\hskip -2em
\endalign
$$
Substituting \thetag{7.1} into \thetag{6.32}, we obtain the determining
equation for the field of point symmetry $U$ in the following form:
$$
F_{\text{det}}(x,y,\ldots,A_{vv},A_{\theta\theta})=0.\hskip -2em
\tag7.2
$$
In left hand side of \thetag{7.2} we have all variables \thetag{6.30}
except for $A_{v\theta}$ The equality \thetag{6.32} brought to the
form \thetag{7.2} should be fulfilled identically in the space
$\Bbb R^{18}$.\par
     The equation \thetag{7.1} is polynomial in derivatives. Therefore
the variables $A_x$, $A_y$, $A_v$, $A_\theta$, $A_{xx}$, $A_{xy}$,
$A_{xv}$, $A_{x\theta}$, $A_{yy}$, $A_{yv}$, $A_{y\theta}$, $A_{vv}$,
$A_{\theta\theta}$ enter the equality \thetag{7.2} polynomially.
The entries of variables $A_{xx}$, $A_{yy}$, $A_{vv}$, and
$A_{\theta\theta}$ in the equality \thetag{7.2} are linear.
In highest order with respect to $A_{xx}$ the equation \thetag{7.2}
looks like
$$
v^3\,A\,\cos\theta\,\bigl(X_\theta+X_A\,A_\theta)\cdot A_{xx}+\ldots=0.
\hskip -2em
\tag7.3
$$
Here $A_\theta$ is independent variable. The quantities $v$, $A$,
$\cos\theta$ in this equation cannot be zero simultaneously. Therefore
from the equality \thetag{7.3} we obtain that the following two partial
derivatives are zero:
$$
\xalignat 2
&\quad X_\theta=\frac{\partial X}{\partial\theta}=0,
&&X_A=\frac{\partial X}{\partial A}=0,\hskip -2em
\tag7.4
\endxalignat
$$
This means that first component $X$ in the field of point symmetry
\thetag{6.3} shouldn't depend on the variables $\theta$ and $A$.
Analyzing the entries of derivatives $A_{yy}$ and $A_{vv}$ in
\thetag{7.2} we obtain analogous results:
$$
\xalignat 2
&\quad Y_\theta=\frac{\partial Y}{\partial\theta}=0,
&&Y_A=\frac{\partial Y}{\partial A}=0,\hskip -2em
\tag7.5\\
\vspace{1ex}
&\quad V_\theta=\frac{\partial V}{\partial\theta}=0,
&&V_A=\frac{\partial V}{\partial A}=0.\hskip -2em
\tag7.6
\endxalignat
$$
Taking into account the relationships \thetag{7.4}, \thetag{7.5}, and
\thetag{7.6}, we can split the equation \thetag{7.2} with respect to
the variable $A_{\theta\theta}$, the entry of which is linear:
$$
F^{\sssize (1)}_{\text{det}}(x,\ldots,A_{y\theta})\cdot A_{\theta
\theta}+F^{\sssize (0)}_{\text{det}}(x,\ldots,A_{y\theta})=0.
$$
From this equation we obtain two separate equations. We write them as
follows:
$$
\aligned
&F^{\sssize (1)}_{\text{det}}(x,y,v,\ldots,A_{yv},A_{y\theta})=0,\\
\vspace{1ex}
&F^{\sssize (0)}_{\text{det}}(x,y,v,\ldots,A_{yv},A_{y\theta})=0.
\endaligned\hskip -2em
\tag7.7
$$
Calculating $F^{\sssize(1)}_{\text{det}}$ explicitly, we find that first
equation \thetag{7.7} doesn't contain actual dependence on $A_{xy}$,
$A_{xv}$, $A_{x\theta}$, $A_{yv}$, $A_{y\theta}$. While with respect to
$A_{\theta}$ left hand side of this equation is a polynomial of the
second order:
$$
-2\,A\,v\,\Theta_A\cdot{A_\theta\vphantom{\tilde A}}^{\!\!2}+\ldots=0.
\hskip -2em
\tag7.8
$$
From the relationship \thetag{7.8} we see that the partial derivative
$\Theta_A$ should vanish:
$$
\Theta_A=\frac{\partial\,\Theta}{\partial A}=0.\hskip -2em
\tag7.9
$$
This means that fourth component of the vector field of point symmetry
\thetag{6.3} do not depend on $A$. With the equality \thetag{7.9} in
mind, the equation \thetag{7.8} turns to be linear in $A_\theta$. We
write it as follows:
$$
F^{\sssize (11)}_{\text{det}}(x,\ldots,A_v)\cdot A_\theta+
F^{\sssize (10)}_{\text{det}}(x,\ldots,A_v)=0.
$$
The above equation splits into two separate equations:
$$
\xalignat 2
&\quad F^{\sssize (11)}_{\text{det}}=0, &&F^{\sssize (10)}_{\text{det}}
=0.\hskip -3em
\tag7.10
\endxalignat
$$
The equations \thetag{7.10} are sufficiently small. We can write them
explicitly. These equations are written as follows:
$$
\align
&\aligned
&V_x\,\cos^2\!\theta+V_y\,\sin^2\!\theta+\frac{A\,V_v}{v}=\\
\vspace{1ex}
&\qquad\quad=\frac{A\,V}{v^2}+\frac{a}{v}-\frac{A\,a_A}{v}+2\,\frac{A\,
\Theta_\theta}{v}.
\endaligned\hskip -3em
\tag7.11\\
\vspace{2ex}
&\Theta_x\,\cos^2\!\theta+\Theta_y\,\sin^2\!\theta+\frac{A\,\Theta_v}{v}
=-a_\theta.\hskip -3em
\tag7.12
\endalign
$$\par
     Now let's consider the second equation \thetag{7.7}. It is linear
with respect to variables $A_{x\theta}$, $A_{y\theta}$. We write it
in the following form:
$$
F^{\sssize (0x)}_{\text{det}}\cdot A_{x\theta}
+F^{\sssize (0y)}_{\text{det}}\cdot A_{y\theta}
+F^{\sssize (00)}_{\text{det}}=0.\hskip -3em
\tag7.13
$$
The equation \thetag{7.13} splits into three separate equations:
$F^{\sssize (0x)}_{\text{det}}=0$, $F^{\sssize (0y)}_{\text{det}}=0$,
and $F^{\sssize (00)}_{\text{det}}=0$. \pagebreak The last equation
doesn't contain second order derivatives $A_{xy}$, $A_{xv}$, $A_{x\theta}$,
$A_{yv}$, $A_{y\theta}$. Left hand side of this equation is a
polynomial of the third order with respect to variable $A_\theta$.
In the leading term we have 
$$
v\,A\,a_{AA}\cdot{A_\theta\vphantom{\tilde A}}^{\!\!3}+\ldots=0.
\hskip -3em
\tag7.14
$$
From the relationship \thetag{7.14} we get $a_{AA}=0$. This means that
fifth component of the vector field of point symmetry \thetag{6.3} is
linear in $A$:
$$
a=\alpha(x,y,v,\theta)+\beta(x,y,v,\theta)\cdot A.\hskip -3em
\tag7.15
$$
Moreover, from the relationships \thetag{7.4}, \thetag{7.5}, \thetag{7.6},
and \thetag{7.9} we get
$$
\xalignat 2
&X=X(x,y,v),&&Y=Y(x,y,v),\hskip -3em\\
\vspace{-1.0ex}
&&&\tag7.17\\
\vspace{-1.5ex}
&V=V(x,y,v),&&\Theta=\Theta(x,y,v,\theta).\hskip -3em
\endxalignat
$$
Thus, the entries of the variable $A$ in components of the field
of point symmetry \thetag{6.3} are calculated in explicit form.\par
     Let's substitute \thetag{7.15} into the equation \thetag{7.11}.
As a result we obtain the equation, which is linear in $A$. It splits
into two separate equations:
$$
\xalignat 2
&\quad V_v=\frac{V}{v}-2\,\Theta_\theta,
&&\alpha=v\,V_x\,\cos\theta+v\,V_y\,\sin\theta.
\hskip -3em
\tag7.18
\endxalignat
$$
Second equation \thetag{7.18} expresses the function $\alpha(x,y,v,
\theta)$ through $V_x$ and $V_y$ in explicit form. Let's substitute
this expression into \thetag{7.15}, then substitute the resulting
expression for $a$ into the equations $F^{\sssize (0x)}_{\text{det}}=0$
and $F^{\sssize (0y)}_{\text{det}}=0$. As a result of such substitution
left hand sides of these equations turn to be polynomials of the second
order with respect to $A$. In leading order they have the following
form:
$$
\xalignat 2
&\quad X_v\cdot A^2+\ldots=0,&&Y_v\cdot A^2+\ldots=0.\hskip -3em
\tag7.19
\endxalignat
$$
From \thetag{7.19} we get $X_v=0$ and $Y_v=0$. Now we can specify
the relationships \thetag{7.17}:
$$
\xalignat 2
&X=X(x,y),&&Y=Y(x,y),\hskip -3em
\tag7.20
\endxalignat
$$
Due to \thetag{7.20} in the equations $F^{\sssize (0x)}_{\text{det}}=0$
and $F^{\sssize (0y)}_{\text{det}}=0$ we can keep only terms linear in
$A$, zero order terms in $A$ being identically zero. First order terms
in these equations can be written explicitly. Thus we can reduce
\thetag{7.19} to the equations that do not contain the variable $A$
at all:
$$
\aligned
&(X_y+\Theta)\,\sin\theta=\left(\!V_v+\frac{V}{v}-\beta-X_x\!\right)
\cos\theta,\\
\vspace{1ex}
&(Y_x-\Theta)\,\cos\theta=\left(\!V_v+\frac{V}{v}-\beta-Y_y\!\right)
\sin\theta.
\endaligned\hskip -3em
\tag7.21
$$
Let's multiply the equation \thetag{7.21} by $\sin\theta$, multiply
second one by $\cos\theta$, then let's subtract second equation from
the first one. At a result we obtain the equation expressing parameter
$\Theta$ through $X_x$, $X_y$, $Y_x$, and $Y_y$:
$$
\Theta=Y_x\,\cos^2\!\theta-X_y\,\sin^2\!\theta+(Y_y-X_x)\,\sin\theta\,
\cos\theta.\hskip -3em
\tag7.22
$$\par
    Further let's consider again the equations \thetag{7.11} and
\thetag{7.12}. First of them is already reduced to the pair of
separate equations \thetag{7.18}, one of which expresses $\alpha$
through $V_x$ and $V_y$. Let's substitute this expression for
$\alpha$ and the expression \thetag{7.22} for $\Theta$ into the
first equation \thetag{7.18}. Upon rather simple transformations
of trigonometric functions we get the equation
$$
(Y_y-X_x)\,\cos(2\theta)-(Y_x+X_y)\,\sin(2\theta)=\frac{V_v}{2}
-\frac{V}{2\,v}.\hskip -3em
\tag7.23
$$
The variable $\theta$ enter to the equation \thetag{7.23} only
through trigonometric functions $\cos(2\theta)$ and $\sin(2\theta)$ 
(see the relationships \thetag{7.17} and \thetag{7.20} above). 
Therefore \thetag{7.23} splits into three separate equations.
One of them is the following
$$
V_v=\frac{V}{v}.\hskip -3em
\tag7.24
$$
Other two equations (it's remarkable) have the form of Cauchy-Riemann
equations with respect to functions $X(x,y)$ and $Y(x,y)$:
$$
\xalignat 2
&X_y=-Y_x,&&X_x=Y_y.\hskip -3em
\tag7.25
\endxalignat
$$
Let's compare \thetag{7.24} and \thetag{7.18}. From such comparison
we get $\Theta_\theta=0$. The equation \thetag{7.24} itself is, in
essential, an ordinary differential equation for $V$. Its general
solution is given by $V=v\,W(x,y)$. Now we can specify \thetag{7.17},
writing these relationships in the following form:
$$
\xalignat 2
&X=X(x,y),&&Y=Y(x,y),\hskip -3em\\
\vspace{-1.0ex}
&&&\tag7.26\\
\vspace{-1.5ex}
&V=v\,W(x,y),&&\Theta=\Theta(x,y,v).\hskip -3em
\endxalignat
$$
Note that the last relationship \thetag{7.26} doesn't contradict
to the formula \thetag{7.22} for $\Theta$. Indeed, if we take
into account Cauchy-Riemann equations \thetag{7.25} and use
trigonometric identity \ $\sin^2\!\theta=1-\cos^2\!\theta$, \ then
we get 
$$
\Theta=-X_y=Y_x.\hskip -3em
\tag7.27
$$
Formula \thetag{7.27} gives further specification of the relationships
\thetag{7.26}:
$$
\xalignat 2
&X=X(x,y),&&Y=Y(x,y),\hskip -3em\\
\vspace{-1.0ex}
&&&\tag7.28\\
\vspace{-1.5ex}
&V=v\,W(x,y),&&\Theta=\Theta(x,y).\hskip -3em
\endxalignat
$$\par
    In the next step we substitute \thetag{7.27} into any one of the
relationships \thetag{7.21}. The results of these two substitutions
coincide with each other. They yield the equality determining
parameter $\beta$ in \thetag{7.15}:
$$
\beta=2\,W-Y_y.\hskip -3em
\tag7.29
$$
From the whole set of the above equations now remain two equations.
These are the equation \thetag{7.12} and the equation
$F^{\sssize (00)}_{\text{det}}=0$, which was derived from \thetag{7.13}.
Let's substitute the reduced functions \thetag{7.28} into the equation
\thetag{7.12}. Thereby we take into account formulas \thetag{7.18} and
\thetag{7.29} for $\alpha$ and $\beta$, formula \thetag{7.27} for
$\Theta$, and formula \thetag{7.15} for parameter $a$. We keep in mind
Cauchy-Riemann equations \thetag{7.25} for $X$ and $Y$ as well. Then
we obtain the equation
$$
(W_y-X_{xy})\,\sin\theta+(X_{xx}-W_x)\,\cos\theta=0.\hskip -3em
\tag7.30
$$
The equation \thetag{7.30} splits into two separate equations, which
form the system of Pfaff equations with respect to the function
$W=W(x,y)$:
$$
\cases W_x=X_{xx},\\W_y=X_{xy}.\endcases\hskip -3em
\tag7.31
$$
System of Pfaff equations \thetag{7.31} is compatible due to the
Cauchy-Riemann equations \thetag{7.25}. It's easy to integrate it.
As a result we determine the function $W(x,y)$ up to some arbitrary
constant $C$:
$$
W=X_x+C.\hskip -3em
\tag7.32
$$
Now it remains to consider only one equation $F^{\sssize(00)}_{\text{det}}
=0$. In a highest order in $A$ it yields $a_{AA}=0$ in \thetag{7.14},
which, in turn, yields \thetag{7.15}. In lower orders in $A$ this equation
appears to be identically fulfilled due to above formulas for the functions
$X$, $Y$, $V$, $\Theta$, and $a$.
\head
\S\,8. Algebra of point symmetries of the reduced normality equation.
\endhead
\rightheadtext{\S\,8. Algebra of point symmetries \dots}
    Results obtained in \S\,7 allow to give complete description of
the structure of any vector field of point symmetry for the equation
\thetag{6.1}. It's convenient to do it in terms of complex variable
$z=x+i\,y$. Due to Cauchy-Riemann equations \thetag{7.25} the components
$X$ and $Y$ in vector field \thetag{6.3} are real and imaginary parts
of some holomorphic function $\xi(z)$. We write them as
$$
\xalignat 2
&\quad X=\Re\xi(z),&&Y=\Im\xi(z).\hskip -2em
\tag8.1
\endxalignat
$$
The relationships \thetag{7.27} and \thetag{7.32} now are written as
follows:
$$
\xalignat 2
&\quad\Theta=\Im\xi'(z),&&W=\Re\xi'(z)+C.\hskip -2em
\tag8.2
\endxalignat
$$
In addition to $z=x+i\,y$ we introduce another complex variable
$u=v\,e^{i\theta}$. Then from the relationships \thetag{7.18}
and \thetag{7.29} we get
$$
\xalignat 2
&\quad\alpha=\Re(u\,\xi''(z)),&&\beta=\Re\xi'(z)+2\,C.\hskip -2em
\tag8.3
\endxalignat
$$
Let's substitute \thetag{8.3} into the formula \thetag{7.15}.
Thereafter we can write formula for vector field of point symmetry
$$
\gather
U=\xi(z)\,\frac{\partial}{\partial z}+\overline{\xi(z)}\,\frac{\partial}
{\partial\bar z}+\xi'(z)\,u\,\frac{\partial}{\partial u}+\overline{\xi'(z)}
\,\bar u\,\frac{\partial}{\partial\bar u}\,+\hskip -2em\\
\vspace{1.4ex}
+\,C\,u\,\frac{\partial}{\partial u}+C\,\bar u\,\frac{\partial}{\partial
\bar u}+|u|\,\frac{\xi''(z)\,u+\overline{\xi''(z)}\,\bar u}{2}\frac{\partial}
{\partial A}\,+\hskip -2em
\tag8.4\\
\vspace{1.4ex}
+\,\frac{\xi'(z)+\overline{\xi'(z)}}{2}\,A\,\frac{\partial}{\partial A}+
2\,C\,A\,\frac{\partial}{\partial A}.\hskip -2em
\endgather
$$
\proclaim{Theorem 8.1} Each vector field $U$ of point symmetry of
reduced normality equation \thetag{6.1} is determined by some holomorphic
function $\xi(z)$ and some real constant $C$ according to the formula
\thetag{8.4}.
\endproclaim
    Let $U_1$ and $U_2$ be two vector fields of point symmetry for
the equation \thetag{6.1}. Suppose that first of them is determined
by the function $\xi_1(z)$ and by constant $C_1$, while second is
determined by function $\xi_2(z)$ and by constant $C_2$. By means of
direct calculations we find that commutator of the fields $U_3=[U_1,
\,U_2]$ has the form \thetag{8.4}. It is determined by holomorphic
function
$$
\xi_3(z)=\xi_1(z)\,\xi_2'(z)-\xi_2(z)\,\xi_1'(z)\hskip -2em
\tag8.5
$$
and by constant $C_3$, where $C_3=0$. We can associate functions $\xi_1(z)$
and $\xi_2(z)$ with the following two holomorphic vector fields:
$$
\xalignat 2
&\quad\xi_1(z)\,\frac{\partial}{\partial z},
&&\xi_2(z)\,\frac{\partial}{\partial z}.\hskip -2em
\tag8.6
\endxalignat
$$
Then function \thetag{8.5} is associated with commutator of holomorphic
vector fields \thetag{8.6}. Adding constants $C_1$ and $C_2$ to vector
fields \thetag{8.6} corresponds to the construction known as
{\it central extension} (see \cite{Brb1}, Chapter~\uppercase
\expandafter{\romannumeral 1}, \S\,1, subsection 7).
\proclaim{Theorem 8.2} Lie algebra $\goth L$ of point symmetries
of normality equation \thetag{6.1}  is isomorphic to central extension
of Lie algebra of holomorphic vector fields in $\Bbb C=\Bbb R^2$. 
\endproclaim
\head
\S\,9. Invariant solutions.
\endhead
    Let $U$ be some vector field of point symmetry of reduced
normality equation \thetag{6.1}. It generates some one-parametric
group of point transformations \thetag{6.4}. Say that function
$A(x,y,v,\theta)$ is {\bf invariant} with respect to vector field
$U$ in \thetag{6.3} if function $A^t(\tilde x,\tilde y,\tilde v,
\tilde\theta)$ defined implicitly by the relationships \thetag{6.10}
do not depend on $t$ and coincides with $A(\tilde x,\tilde y,
\tilde v,\tilde\theta)$. This condition is expressed by the
relationship
$$
\tilde A=A(\tilde x,\tilde y,\tilde v,\tilde\theta),\hskip -2em
\tag9.1
$$
which binds functions $\tilde x$, $\tilde y$, $\tilde v$, $\tilde\theta$,
and $\tilde A$ in \thetag{6.10}. Let's differentiate the relationship
\thetag{9.1} with respect to parameter $t$. Then substitute $t=0$ and
take into account the expansions \thetag{6.7}. As a result we get the
relationship 
$$
\gather
A_x(x,y,v,\theta)\cdot X(x,y,v,\theta,A(x,y,v,\theta))\,+\\
\ \hskip 2cm+\,A_y(x,y,v,\theta)\cdot Y(x,y,v,\theta,A(x,y,v,\theta))\,+\\
+\,A_v(x,y,v,\theta)\cdot V(x,y,v,\theta,A(x,y,v,\theta))\,+\\
\ \hskip 2cm+\,A_\theta(x,y,v,\theta)\cdot\Theta(x,y,v,\theta,
A(x,y,v,\theta))=\\
\ \hskip 1cm=a(x,y,v,\theta,A(x,y,v,\theta)).
\endgather
$$
This relationship expresses the condition of {\bf invariance} of the
function $A(x,y,v,\theta)$ with respect to vector field \thetag{6.3}
in the form that can be checked by direct calculations. Briefly
this relationship is written as 
$$
A_x\cdot X+A_y\cdot Y+A_v\cdot V+A_\theta\cdot\Theta=a.\hskip -2em
\tag9.2
$$
\definition{Definition 9.1} The solution $A(x,y,v,\theta)$ of reduced
normality equation \thetag{6.1} is called {\bf invariant solution}
if it is invariant with respect to one or several fields of point symmetry
of this equation.
\enddefinition
    Vector fields of point symmetry of the equation \thetag{6.1}, with
respect to which some fixed solution $A(x,y,v,\theta)$ of this
equation is invariant, form some Lie subalgebra $\goth L_A$ in Lie
algebra $\goth L$ of all point symmetries of this equation. All examples
of dynamical systems admitting the normal shift, which were considered
in \S\,4, correspond to some invariant solutions of reduced normality
equation. Now we consider these examples again and calculate Lie
subalgebras corresponding to them.
\demo{Example 1} {\bf Spatially homogeneous, but not isotropic force
field}. Function $A(x,y,v,\theta)$ in this example has the form
$A=A(v)\,\cos\theta$. Let's substitute this function into the equation
\thetag{9.2}. Thereby we take into account the following formulas for 
the components of the field of point symmetry:
$$
\align
&V=(X_x+C)\,v,\hskip -2em\\
\vspace{1ex}
&\Theta=-X_y,\hskip -2em
\tag9.3\\
\vspace{1ex}
&a=v^2\,(X_{xx}\,\cos\theta+X_{xy}\,\sin\theta)+(X_x+2\,C)\,A.
\hskip -2em
\endalign
$$
These formulas are obtained from the relationships \pagebreak \thetag{7.15},
\thetag{7.18}, \thetag{7.27}, \thetag{7.28}, \thetag{7.29}, and
\thetag{7.32}, which were derived in \S\,7. Substituting \thetag{9.3}
and $A=A(v)\,\cos\theta$ into the equation \thetag{9.2}, we obtain the
equation:
$$
\gathered
v^2\,(X_{xx}\,\cos\theta+X_{xy}\,\sin\theta)
+(X_x+2\,C)\,A(v)\cos\theta=\\
\vspace{1ex}
=v\,A'(v)\,(X_x+C)\,\cos\theta+A(v)\,X_y\,\sin\theta.
\endgathered\hskip -2em
\tag9.4
$$
By collecting coefficients of $\cos\theta$ and $\sin\theta$
in \thetag{9.4} we split this equation into two separate equations.
They are written as follows:
$$
\align
&X_{xy}\,v^2=A(v)\,X_y,\hskip -2em
\tag9.5\\
\vspace{1ex}
&X_{xx}\,v^2+(X_x+2\,C)\,A(v)=v\,A'(v)\,(X_x+C).
\hskip -2em
\tag9.6
\endalign
$$
Considering the equations \thetag{9.5} and \thetag{9.6} we find
three subcases.\par
{\bf Subcase 1}. Function $A(v)$ is not an exponential function,
i\.~e\. $A(v)\neq k\cdot v^\alpha$, where $k$ and $\alpha$ are
constants. In this case from \thetag{9.5} we get $X_y=0$. Then
from Cauchy-Riemann equations \thetag{7.25} we derive $Y_x=0$
and $X_x=Y_y=C_1=\const$. Hence for $X$ and $Y$ we have the
expressions
$$
\xalignat 2
&\quad X=C_1\,x+C_2,&&Y=C_1\,y+C_3.\hskip -2em
\tag9.7
\endxalignat
$$
Substituting \thetag{9.7} into \thetag{9.6} and taking into
account that $A(v)\neq c\cdot v^\alpha$, we find that $C=C_1=0$.
Constants $C_2$ and $C_3$ remain undetermined. This means that
Lie subalgebra $\goth L_A$ for the solution $A=A(v)\,\cos\theta$
in this case is two-dimensional. It is generated by the following
two vector fields:
$$
\xalignat 2
&U_1=\frac{\partial}{\partial x},
&&U_2=\frac{\partial}{\partial y}.\hskip -2em
\tag9.8
\endxalignat
$$\par
{\bf Subcase 2}. Function $A(v)$ is exponential function $A(v)=k\cdot
v^\alpha$, where $k$ and $\alpha$ are constants, but $\alpha\neq 2$.
In this case, as in previous one, the equation \thetag{9.4} splits into
two separate equations \thetag{9.5} and \thetag{9.6}. From \thetag{9.5}
we derive \thetag{9.7}. However, the result of substituting \thetag{9.7}
into \thetag{9.6} is different:
$$
C_1\,(1-\alpha)=C\,(\alpha-2).\hskip -2em
\tag9.9
$$
The equation \thetag{9.9} is solvable in form of the following two
relationships:
$$
\xalignat 2
&\quad C_1=(\alpha-2)\,C_4,&&C=(1-\alpha)\,C_4.\hskip -3em
\tag9.10
\endxalignat
$$
Formulas \thetag{9.7} and \thetag{9.10} contain three arbitrary
constants $C_2$, $C_3$, $C_4$. \pagebreak
This indicate that algebra $\goth L_A$
is three-dimensional. Here are the generators of this algebra:
$$
\align
&U_1=\frac{\partial}{\partial x},\qquad\qquad
U_2=\frac{\partial}{\partial y},\hskip -3em
\vspace{-1.0ex}
&\tag9.11\\
\vspace{-1.5ex}
&U_3=(\alpha-2)\left(\!x\,\frac{\partial}{\partial x}+y\,
\frac{\partial}{\partial y}\right)
-v\,\frac{\partial}{\partial v}-\alpha\,k\,v^\alpha\,\cos\theta\,
\frac{\partial}{\partial A}.\hskip -3em
\endalign
$$\par
{\bf Subcase 3}. Function $A(v)$ is quadratic, i\.~e\. $A(v)=k\cdot v^2$. 
The equations \thetag{9.5} and \thetag{9.6} in this case are written as
follows:
$$
\xalignat 2
&X_{xy}=k\,X_y,&&X_{xx}=k\,X_x.\hskip -3em
\tag9.12
\endxalignat
$$
The equations \thetag{9.12} should be completed by the equations of
Cauchy-Riemann \thetag{7.25}. Upon rather simple analysis we can write
general solution to the system of equations for $X$ and $Y$ that arises
in this case:
$$
\aligned
&X=C_1\,\cos(ky)\,e^{kx}+C_2\,\sin(ky)\,e^{kx}+C_3,\\
\vspace{1ex}
&Y=C_2\,\cos(ky)\,e^{kx}-C_1\,\sin(ky)\,e^{kx}+C_4.
\endaligned\hskip -3em
\tag9.13
$$
Elements of algebra $\goth L_A$ depend on five arbitrary constants
$C_1$, $C_2$, $C_3$, $C_4$, and $C$. Here are the generators of
this algebra:
$$
\align
&U_1=\frac{\partial}{\partial x},\qquad\qquad
U_2=\frac{\partial}{\partial y},\hskip -3em\\
\vspace{2ex}
&U_3=v\,\frac{\partial}{\partial v}-2\,k\,v^2\,\cos\theta
\frac{\partial}{\partial A},\hskip -3em\\
\vspace{2ex}
&U_4=e^{kx}\,\cos(ky)\!\left(\frac{\partial}{\partial x}+
k\,v\,\frac{\partial}{\partial v}+2\,k^2\,v^2\,\cos\theta\,
\frac{\partial}{\partial A}\right)+\hskip -3em
\tag9.14\\
\vspace{1ex}
&\qquad\quad+\,e^{kx}\,\sin(ky)\!\left(\frac{\partial}{\partial y}
+k\,\frac{\partial}{\partial\theta}-k^2\,v^2\,\sin\theta\,
\frac{\partial}{\partial A}\right),\hskip -3em\\
\vspace{2ex}
&U_5=e^{kx}\,\sin(ky)\!\left(\frac{\partial}{\partial x}+
k\,v\,\frac{\partial}{\partial v}+2\,k^2\,v^2\,\cos\theta\,
\frac{\partial}{\partial A}\right)-\hskip -3em\\
\vspace{1ex}
&\qquad\quad-\,e^{kx}\,\cos(ky)\!\left(\frac{\partial}{\partial y}
+k\,\frac{\partial}{\partial\theta}-k^2\,v^2\,\sin\theta\,
\frac{\partial}{\partial A}\right).\hskip -3em
\endalign
$$
\enddemo
\demo{Example 2} {\bf Spatially non-homogeneous force field with
marked point}. Function $A$ determining force field $\bold F$ 
in this case is given by formula \thetag{4.29}:
$$
A=\frac{A(v)\,\cos\theta}{\rho}.\hskip -3em
\tag9.15
$$
However, this function is written in variables $\rho$, $\gamma$, $v$,
$\theta$. They are different from variables $x$, $y$, $v$, $\theta$,
in which the reduced normality equation \thetag{6.1} is written.
One easily find the relations between these two sets of variables
by comparing Fig\.~4.1 and Fig\.~4.2. The variable $v$ denotes the
modulus of velocity in both sets of variables. Angular variables
differ from each other by the value of the angle $\gamma$. Therefore
function \thetag{9.15} should be rewritten as follows:
$$
A=\frac{A(v)\,\cos(\theta-\gamma)}{\rho}.
$$
In order to complete the transfer to the variables $x$, $y$, $v$, $\theta$
one should apply the identity \ $\cos(\theta-\gamma)=\cos\theta\,\cos\gamma
+\sin\theta\,\sin\gamma$ \ and one should use the relationships binding
$\rho$ and $\gamma$ with variables $x$ and $y$:
$$
\xalignat 2
&x=\rho\,\cos\gamma,&&y=\rho\,\sin\gamma.
\endxalignat
$$
As a result one get the following expression for $A$:
$$
A=A(v)\,\frac{x\,\cos\theta+y\,\sin\theta}{x^2+y^2}.\hskip -2em
\tag9.16
$$
Let's substitute \thetag{9.16} into the equation \thetag{9.2} and
let's take into account the relationships \thetag{9.3}. Instead of
\thetag{9.5} and \thetag{9.6} we get two other equations
$$
\align
&\aligned
-v^2\,X_{xx}+A(v)&\left(\frac{X-x\,(X_x+2\,C)-y\,X_y}{\rho}\right.-\\
\vspace{1ex}
&-\left.2\,x\,\frac{x\,X+y\,Y}{\rho^2}\right)
+v\,A'(v)\,x\,\frac{X_x+C}{\rho}=0,
\endaligned\hskip -3em
\tag9.17\\
\vspace{2ex}
&\aligned
-v^2\,X_{xy}+A(v)&\left(\frac{Y-y\,(X_x+2\,C)-x\,X_y}{\rho}\right.-\\
\vspace{1ex}
&-\left.2\,y\,\frac{x\,X+y\,Y}{\rho^2}\right)
+v\,A'(v)\,y\,\frac{X_x+C}{\rho}=0.
\endaligned\hskip -3em
\tag9.18
\endalign
$$
Here $\rho=\sqrt{x^2+y^2}$. The equations \thetag{9.17}
and \thetag{9.18} determine three subcases. \par
{\bf Subcase 1}. Function $A(v)$ is not an exponential function, i\.~e\.
$A(v)\neq k\cdot v^\alpha$, where $k$ and $\alpha$ are constants.
In this case each of the above two equation \thetag{9.17} and \thetag{9.18}
splits into three separate equations. Among six resulting equations we have
two coincident ones. So it remains five equations that follows from
\thetag{9.17} and \thetag{9.18}. Here are three of these five equations:
$$
\pagebreak
\xalignat 3
&\quad X_{xx}=0,&&X_{xy}=0,&&X_x=-C.\hskip -3em
\tag9.19
\endxalignat
$$
Other two equations are more complicated. They have the following form:
$$
\aligned
&(X-x\,C-y\,X_y)(x^2+y^2)=2\,x\,(x\,X+y\,Y),\\
\vspace{1ex}
&(Y-y\,C+x\,X_y)(x^2+y^2)=2\,y\,(x\,X+y\,Y).
\endaligned\hskip -3em
\tag9.20
$$
Let's complete \thetag{9.19} with the Cauchy-Riemann equations 
\thetag{7.25}. For $X$ and $Y$, upon solving these equations, we
get
$$
\xalignat 2
&X=-C\,x+C_1\,y+C_2,
&&Y=-C\,y-C_1\,x+C_3.
\endxalignat
$$
Substituting these expressions into \thetag{9.20} we get the equations
that yield $C_2=0$ and $C_3=0$. Constants $C$ and $C_1$ remain
undetermined. Hence Lie subalgebra $\goth L_A$ in this case is
two-dimensional. Here are its generators:
$$
\aligned 
&U_1=x\,\frac{\partial}{\partial x}+y\,\frac{\partial}{\partial y}-
A(v)\,\frac{x\,\cos\theta+y\,\sin\theta}{x^2+y^2}\,\frac{\partial}
{\partial A},\\
\vspace{1ex}
&U_2=-y\,\frac{\partial}{\partial x}+x\,\frac{\partial}{\partial y}
+\frac{\partial}{\partial\theta}.
\endaligned\hskip -3em
\tag9.21
$$\par
{\bf Subcase 2}. Function $A(v)$ is not an exponential function, i\.~e\.
$A(v)\neq k\cdot v^\alpha$, where $k$ and $\alpha$ are constants, but
$\alpha\neq 2$. From \thetag{9.17} and \thetag{9.18} in this case we
extract four equations. Here are two of them:
$$
\xalignat 2
&X_{xx}=0,&&X_{xy}=0.\hskip -3em
\tag9.22
\endxalignat
$$
Rest two equations are written in the following form:
$$
\align
&X+x\,\alpha\,(X_x+C)-x\,(X_x+2\,C)
=y\,X_y+\frac{2\,x\,(x\,X+y\,Y)}{x^2+y^2},
\hskip -3em
\tag9.23\\
\vspace{2ex}
&Y+y\,\alpha\,(X_x+C)-y\,(X_x+2\,C)
=-x\,X_y+\frac{2\,y\,(x\,X+y\,Y)}{x^2+y^2}.
\hskip -3em
\tag9.24
\endalign
$$
Let's complete the equations \thetag{9.22} with Cauchy-Riemann
equations \thetag{7.25}. For $X$ and $Y$ this yields the
following expressions:
$$
\xalignat 2
&X=C_1\,x+C_2\,y+C_3,
&&Y=C_1\,y-C_2\,x+C_4.
\endxalignat
$$
Substituting these expressions into \thetag{9.23} and \thetag{9.24},
we obtain three equations, solution of which is given by the
relationships
$$
\pagebreak
\xalignat 3
&C=-C_1,&&C_3=0,&&C_4=0.
\endxalignat
$$
Constants $C_1$ and $C_2$ remain undetermined, this means that
Lie algebra $\goth L_A$ is two-dimensional: $\dim\goth L_A=2$.
Here are its generators:
$$
\aligned 
&U_1=x\,\frac{\partial}{\partial x}+y\,\frac{\partial}{\partial y}-
k\,v^\alpha\,\frac{x\,\cos\theta+y\,\sin\theta}{x^2+y^2}\,\frac{\partial}
{\partial A},\\
\vspace{1ex}
&U_2=-y\,\frac{\partial}{\partial x}+x\,\frac{\partial}{\partial y}
+\frac{\partial}{\partial\theta}.
\endaligned\hskip -3em
\tag9.25
$$
They have the same form as generators \thetag{9.21} in previous case.
\par
{\bf Subcase 3}. Function $A(v)$ is quadratic, i\.~e\. $A(v)=k\cdot v^2$. 
The equations \thetag{9.17} and \thetag{9.18} in this case do not split.
In order to write these equations and solve them we use complex form
\thetag{8.4} for the vector field of point symmetry \thetag{6.3}.
Let's multiply \thetag{9.18} by $i$ and let's subtract the result from
\thetag{9.17}. The resulting equation leads to the following differential
equation for holomorphic function $\xi(z)$ in \thetag{8.4}:
$$
\xi''(z)-\frac{k}{z}\,\xi'(z)+\frac{k}{z^2}\,\xi(z)=0.\hskip -3em
\tag9.26
$$
The space of solution of linear ordinary differential equation
\thetag{9.26} over complex numbers is two-dimensional:
$$
\xi(z)=C_1\,z^k+C_2\,z.\hskip -3em
\tag9.27
$$
Its dimension over real numbers is equal to $4$. Besides $\xi(z)$
formula \thetag{8.4} contain real constant $C$, the value of which
is not restricted by the equation \thetag{9.26}. Therefore the
dimension of algebra $\goth L_A$ is equal to $5$:
$$
\pagebreak
\aligned
&U_1=x\,\frac{\partial}{\partial x}+y\,\frac{\partial}{\partial y}-
k\,v^2\,\frac{x\,\cos\theta+y\,\sin\theta}{x^2+y^2}\,\frac{\partial}
{\partial A},\\
\vspace{1ex}
&U_2=-y\,\frac{\partial}{\partial x}+x\,\frac{\partial}{\partial y}
+\frac{\partial}{\partial\theta},\\
\vspace{1ex}
&U_3=x\,\frac{\partial}{\partial x}+y\,\frac{\partial}{\partial y}
+\frac{1}{2}\,v\,\frac{\partial}{\partial v},\\
\vspace{1ex}
&\aligned U_4&=X_k\,\frac{\partial}{\partial x}+Y_k\,\frac{\partial}
{\partial y}+k\,v\,X_{k-1}\,\frac{\partial}{\partial v}+k\,Y_{k-1}
\,\frac{\partial}{\partial\theta}\,+\\
&+\,v^2\,k\left((k-1)\,(X_{k-2}\,\cos\theta+Y_{k-2}\,\sin\theta)
\vphantom{\frac{x\,cos\theta+y\,\sin\theta}{x^2+y^2}}
+\frac{x\,cos\theta+y\,\sin\theta}{x^2+y^2}\,k\,
X_{k-1}\right)\frac{\partial}{\partial A},
\endaligned\\
\vspace{1ex}
&\aligned U_5&=-Y_k\,\frac{\partial}{\partial x}+X_k\,\frac{\partial}
{\partial y}-k\,v\,Y_{k-1}\,\frac{\partial}{\partial v}+k\,X_{k-1}
\,\frac{\partial}{\partial\theta}\,+\\
&+\,v^2\,k\left((1-k)\,(Y_{k-2}\,\cos\theta+X_{k-2}\,\sin\theta)
\vphantom{\frac{x\,cos\theta+y\,\sin\theta}{x^2+y^2}}
-\frac{x\,cos\theta+y\,\sin\theta}{x^2+y^2}\,k\,
Y_{k-1}\right)\frac{\partial}{\partial A}.
\endaligned
\endaligned
$$
In writing components of vector fields $U_4$ and $U_5$ we used the
following notations:
$$
\xalignat 2
&X_q(x,y)=\Re\left((x+i\,y)^q\right),
&&Y_q(x,y)=\Im\left((x+i\,y)^q\right).
\endxalignat
$$
\enddemo
\head
\S\,10. Symmetry reduction
of reduced normality equation.
\endhead
    By means of scalar ansatz \thetag{2.4} the system of normality
equations \thetag{1.1} were reduced to one differential equation.
In variable $x$, $y$, $v$, $\theta$ this equation looks like
\thetag{6.1}. Further reduction of this equation in general case
is impossible. However, there are special cases, when such reduction
is possible. Below we consider some of these special cases. 
\demo{Example 3} {\bf Spatially homogeneous, but not isotropic force
field}. For symmetry reduction of normality equation \thetag{6.1}
we need to use some subalgebra $\goth L_A$ in the algebra of point
symmetries $\goth L$ of this equation. In this example we choose
subalgebra $\goth L_A$ with generators \thetag{9.8}:
$$
\xalignat 2
&U_1=\frac{\partial}{\partial x},
&&U_2=\frac{\partial}{\partial y}.\hskip -3em
\tag10.1
\endxalignat
$$
Substituting components of vector fields \thetag{10.1} into the
equation \thetag{9.2}, we get $A_x=0$ and $A_y=0$. This reduces
normality equation \thetag{6.1} to the form \thetag{4.6}:
$$
\frac{A\,A_\theta}{v^2}+\frac{A_\theta\,A_{\theta\theta}}{v^2}
+\frac{A_\theta\,A_v}{v}=\frac{A\,A_{\theta v}}{v}.\hskip -3em
\tag10.2
$$
We can lower the order of differential equation \thetag{10.2}.
Let's denote
$$
b=\frac{A_\theta}{A}.\hskip -3em
\tag10.3
$$
Substituting \thetag{10.3} into \thetag{10.2}, we reduce \thetag{10.2}
to the following quasilinear differential equation of the first order
for the function $b$:
$$
b\,b_\theta-v\,b_v+b^3+b=0.\hskip -3em
\tag10.4
$$
The equation \thetag{10.4} is solved by method of characteristics
(see \cite{Kar1}). Characteristic lines of the equation \thetag{10.4}
are parametric lines in $\Bbb R^3$ given by the equations 
$$
\xalignat 3
&\quad\dot v=-v,&&\dot\theta=b,&&\dot b=-b^3-b.\hskip -3em
\tag10.5
\endxalignat
$$
It's easy to write two functionally independent first integrals
for the system of differential equations \thetag{10.5}. These are 
$$
\xalignat 2
&I_1=\theta+\arctan(b),
&&I_2=\frac{u\,b}{v\sqrt{1+b^2}},\hskip -3em
\tag10.6
\endxalignat
$$
where $u$ is some constant. \pagebreak According to the theory of
quasilinear partial differential equations of the first order
(see \cite{Kar1}),
general solution of the equation \thetag{10.4} is given in implicit
form by functional equation 
$$
\varPhi(I_1,I_2)=0,\hskip -3em
\tag10.7
$$
where $\varPhi$ is some arbitrary function of two variables.
From \thetag{10.6} and \thetag{10.7} it follows that the equation
\thetag{10.4} is integrable in quadratures. Hence the equation
\thetag{10.2}  is also integrable in quadratures. This result
was obtained in \cite{Bol8}.\enddemo
\demo{Example 4} {\bf Spatially non-homogeneous force field with
marked point}. In this example we choose Abelian
Lie algebra $\goth L_A$ with generators \thetag{9.21}:
$$
\xalignat 2
&U_1=x\,\frac{\partial}{\partial x}+y\,\frac{\partial}{\partial y}
-A\,\frac{\partial}{\partial A},
&&U_2=-y\,\frac{\partial}{\partial x}+x\,\frac{\partial}{\partial y}
+\frac{\partial}{\partial\theta}.
\endxalignat
$$
Substituting components of these two vector fields into the
equation \thetag{9.2}, we get the system of two differential
equations:
$$
\aligned
&A_x-\frac{y}{x^2+y^2}\,A_\theta+\frac{x}{x^2+y^2}\,A=0,\\
\vspace{1ex}
&A_y+\frac{x}{x^2+y^2}\,A_\theta+\frac{y}{x^2+y^2}\,A=0.
\endaligned\hskip -3em
\tag10.9
$$
The equations \thetag{10.9} are compatible. Similar to \thetag{10.4},
the equations \thetag{10.9} are quasilinear equations of the first order.
Such equations are solved by method of characteristics (see \cite{Kar1}).
Applying this method to \thetag{10.9}, we get the theorem.
\proclaim{Theorem 10.1} Solution of the system of equations \thetag{10.4}
is determined by an arbitrary function of two variables $A(v,\theta)$:
$$
A=\frac{A(v,\theta-\gamma)}{\rho}.\hskip -3em
\tag10.10
$$
Here $\rho=\sqrt{x^2+y^2}$, while the angle $\gamma$ in formula
\thetag{10.10} is determined by the relationships \ $x=\rho\,cos\,\gamma$
\ and \ $y=\rho\,sin\,\gamma$.
\endproclaim
    Relying upon theorem~10.1, we substitute \thetag{10.10} into the
normality equation \thetag{6.1}. This determines symmetry reduction
of the equation \thetag{6.1}. In this case it's given by the
following equation for the function $A=A(v,\theta)$:
$$
\frac{A\,A_\theta}{v^2}+\frac{A_\theta\,A_{\theta\theta}}{v^2}
+\frac{A_\theta\,A_v}{v}+(A+A_{\theta\theta})\,\sin\theta=
\frac{A\,A_{\theta v}}{v}.\hskip -1em
\tag10.11
$$
The equation \thetag{10.11} is quite similar to \thetag{10.2}. But
the substitution \thetag{10.3} doesn't lower the order of the equation
\thetag{10.11}.\enddemo
\demo{Example 5} {\bf Reduction determined by algebra of conformal
automorphisms of a disc}. In this example we consider one more way
for reducing the normality equation \thetag{6.1}. Here we use
complex form of vector field of point symmetry \thetag{8.4} and
theorem~8.2, which binds algebra $\goth L$ with algebra of holomorphic
vector fields in $\Bbb C=\Bbb R^2$. Let's consider the disc $K$ of
the radius $R$ in complex plane:
$$
K=\{z\in\Bbb C: \ |z|<R\}.
$$
Group of conformal automorphisms of the domain $K$ has the real
dimension $3$. Each automorphism in this group is determined by
one complex parameter $\alpha$ and one real parameter $\varphi$
(see more details in \cite{Shb1}):
$$
f(z)=e^{i\varphi}\,\frac{R\,(z-R\,\alpha)}{R-\bar\alpha\,z}.\hskip -3em
\tag10.12
$$
For $\varphi\to 0$ and $\alpha\to 0$ we get $f(z)\to\id(z)$. Let's
calculate the derivatives of the function $f(z)$ with respect to
parameters $\varphi$, $\alpha$, and $\bar\alpha$ at the point
$\varphi=0$ and $\alpha=0$:
$$
\gather
\xi_1(z)=\frac{\partial f(z)}{\partial\varphi}\,\hbox{\vrule height
12pt depth 8pt width 0.5pt}_{\,\aligned &\ssize\varphi=0\\ \vspace{-2ex}
&\ssize\alpha=0\endaligned}\ =\ i\,z,\\
\vspace{1ex}
\xi_2(z)=\frac{\partial f(z)}{\partial\alpha}\,
\hbox{\vrule height 12pt depth 8pt width 0.5pt}_{\,\aligned
&\ssize\varphi=0\\\vspace{-2ex}&\ssize\alpha=0\endaligned}\ +\
\frac{\partial f(z)}{\partial\bar\alpha}\,\hbox{
\vrule height 12pt depth 8pt width 0.5pt}_{\,\aligned &\ssize
\varphi=0\\\vspace{-2ex}&\ssize\alpha=0\endaligned}\ =
\frac{z^2-R^2}{R},\\
\vspace{1ex}
\xi_3(z)=i\,\frac{\partial f(z)}{\partial\alpha}\,
\hbox{\vrule height 12pt depth 8pt width 0.5pt}_{\,\aligned
&\ssize\varphi=0\\\vspace{-2ex}&\ssize\alpha=0\endaligned}\ -\
i\,\frac{\partial f(z)}{\partial\bar\alpha}\,\hbox{
\vrule height 12pt depth 8pt width 0.5pt}_{\,\aligned&\ssize
\varphi=0\\\vspace{-2ex}&\ssize\alpha=0\endaligned}\ =-i\,
\frac{z^2+R^2}{R}.
\endgather
$$
Functions $\xi_1$, $\xi_2$, and $\xi_3$ determine three holomorphic
vector fields, which generate algebra of conformal automorphism of
a disk $K$:
$$
\xalignat 3
&\xi_1(z)\,\frac{\partial}{\partial z},
&&\xi_2(z)\,\frac{\partial}{\partial z},
&&\xi_3(z)\,\frac{\partial}{\partial z}.
\endxalignat
$$
It's easy to calculate commutators of these vector fields:
$$
\xalignat 3
&[\xi_1,\,\xi_2]=-\xi_3,
&&[\xi_2,\,\xi_3]=4\,\xi_1,
&&[\xi_3,\,\xi_1]=-\xi_2.
\endxalignat
$$
Let's substitute the functions $\xi_1(z)$, $\xi_2(z)$, and $\xi_3(z)$
into the formula \thetag{8.4} for vector field of point symmetry of
the equation \thetag{6.1}, taking $C=0$ by each such substitution.
This determines three fields $U_1$, $U_2$, $U_3$, which are generators
of three-dimensional non-abelian subalgebra $\goth L_A$ in the algebra
of point symmetries of the equation \thetag{6.1}:
$$
\xalignat 3
&[U_1,\,U_2]=-U_3,
&&[U_2,\,U_3]=4\,U_1,
&&[U_3,\,U_1]=-U_2.
\endxalignat
$$
In order to write explicit formulas for vector fields $U_1$, $U_2$, $U_3$
in real variables $x$, $y$, $v$, $\theta$ we use formulas \thetag{8.1},
\thetag{8.2}, \thetag{8.2}, and formulas \thetag{7.15}, \thetag{7.26}:
$$
\align
&U_1=-y\,\frac{\partial}{\partial x}+x\,\frac{\partial}{\partial y}
+\frac{\partial}{\partial\theta},\hskip -3em 
\tag10.13\\
\vspace{2ex}
&\aligned
 U_2&=\frac{x^2-y^2-R^2}{R}\,\frac{\partial}{\partial x}+\frac{2\,x\,y}
 {R}\,\frac{\partial}{\partial y}\,+\\
 \vspace{1ex}
 &+\,\frac{2\,v\,x}{R}\,\frac{\partial}{\partial v}+\frac{2\,y}{R}
 \,\frac{\partial}{\partial\theta}+2\,\frac{v^2\,\cos\theta+x\,A}
 {R}\,\frac{\partial}{\partial A},
\endaligned\hskip -3em
\tag10.14\\
\vspace{2ex}
&\aligned
 U_3&=\frac{2\,x\,y}{R}\,\frac{\partial}{\partial x}-\frac{x^2-y^2+R^2}
 {R}\,\frac{\partial}{\partial y}\,+\\
 \vspace{1ex}
 &+\,\frac{2\,v\,y}{R}\,\frac{\partial}{\partial v}-\frac{2\,x}{R}\,
 \frac{\partial}{\partial\theta}+2\,\frac{v^2\,\sin\theta+y\,A}{R}\,
 \frac{\partial}{\partial A}.
\endaligned\hskip -3em
\tag10.15
\endalign
$$
Lie algebra $\goth L_A$ with generators \thetag{10.13}, \thetag{10.14},
\thetag{10.15} is isomorphic to matrix algebra $\so(1,2,\Bbb R)$.\par
     Vector field $U_1$ in \thetag{10.13} coincides with the field
$U_2$ in the above example~2 (see formulas \thetag{9.21}). For to
simplify formulas \thetag{10.13}, \thetag{10.14}, \thetag{10.15}
we transfer from $x$, $y$, $v$, $\theta$ to variables $\rho$, $\gamma$,
$v$, $\theta$ used in the example~2. In these variables vector field
\thetag{10.13} is reduced to the operator of differentiation in
$\gamma$:
$$
U_1=\frac{\partial}{\partial\gamma}.\hskip -3em
\tag10.16
$$
While normality equation \thetag{6.1} turns to \thetag{4.28}. The
invariance condition of the solution of \thetag{4.28} with respect
to the field \thetag{10.16} is written as $A_\gamma=0$. This means
that $A=A(\rho,v,\theta)$. Let's transform vector fields \thetag{10.14}
and \thetag{10.15} to the variables $\rho$, $\gamma$, $v$, $\theta$.
For vector field $U_2$ we get
$$
\aligned
U_2&=\frac{\rho^2-R^2}{R}\left(\cos\gamma\,\frac{\partial}{\partial\rho}+
\frac{\sin\gamma}{\rho}\,\frac{\partial}{\partial\theta}\right)+
\frac{2\,v\,\rho\cos\gamma}{R}\,\frac{\partial}{\partial v}\,+\\
\vspace{1ex}
&+\,\frac{\rho^2+R^2}{R}\,\frac{\sin\gamma}{\rho}\,\frac{\partial}
{\partial\gamma}+2\,\frac{\cos(\theta+\gamma)\,v^2+A\,\rho\,\cos\gamma}{R}
\,\frac{\partial}{\partial A}.
\endaligned
$$
For the field $U_3$ we get similar expression: 
$$
\pagebreak
\aligned
U_3&=\frac{\rho^2-R^2}{R}\left(\sin\gamma\,\frac{\partial}{\partial\rho}-
\frac{\cos\gamma}{\rho}\,\frac{\partial}{\partial\theta}\right)+
\frac{2\,v\,\rho\sin\gamma}{R}\,\frac{\partial}{\partial v}\,-\\
\vspace{1ex}
&-\,\frac{\rho^2+R^2}{R}\,\frac{\cos\gamma}{\rho}\,\frac{\partial}
{\partial\gamma}+2\,\frac{\sin(\theta+\gamma)\,v^2+A\,\rho\,\sin\gamma}{R}
\,\frac{\partial}{\partial A}.
\endaligned
$$
Writing for $A=A(\rho,v,\theta)$ the invariance conditions with respect
to $U_2$ and $U_3$, we can reduce these conditions to the form of
equations
$$
\align
&A_\theta=\frac{2\,\rho\,v^2\,\sin\theta}{R^2-\rho^2},\hskip -3em
\tag10.17\\
\vspace{2ex}
&A_\rho=\frac{2\,\rho\,v\,A_v-2\,v^2\,\cos\theta-2\,\rho\,A}{R^2-\rho^2}.
\hskip -3em
\tag10.18
\endalign
$$
The equation \thetag{10.17} is easily solved. This determine the
dependence of $A$ on $\theta$:
$$
A=-\frac{2\,\rho\,v^2\,\cos\theta}{R^2-\rho^2}+A(\rho,v).\hskip -3em
\tag10.19
$$
Let's substitute \thetag{10.19} into \thetag{10.18}. As a result we get
the following equation for the function $A(\rho,v)$ from \thetag{10.19}:
$$
A_\rho-\frac{2\,\rho\,v}{R^2-\rho^2}\,A_v+\frac{2\,\rho\,A}{R^2-\rho^2}
=0.\hskip -3em
\tag10.20
$$
The equation is solved by method of characteristics (see \cite{Kar1}).
General solution of this equation contain one arbitrary function $f=f(v)$:
$$
A(\rho,v)=v\,f\!\left(\frac{v}{R^2-\rho^2}\right).
$$
Now let's substitute the above function into \thetag{10.19}. Invariant
solution of the equation \thetag{4.28} then should be found in the
following form:
$$
A=-\frac{2\,\rho\,v^2\,\cos\theta}{R^2-\rho^2}+
v\,f\!\left(\frac{v}{R^2-\rho^2}\right).\hskip -3em
\tag10.21
$$
Substituting \thetag{10.21} into \thetag{4.28}, we find that
\thetag{10.21} is the solution of the equation \thetag{4.28}
for any choice of function $f=f(v)$. In variables $x$, $y$, $v$,
$\theta$ for $A$ we have
$$
A=-\frac{2\,v^2\,(x\,\cos\theta+y\,\sin\theta)}{R^2-x^2-y^2}+
v\,f\!\left(\frac{v}{R^2-x^2-y^2}\right).\hskip -3em
\tag10.22
$$
Function \thetag{10.22} is the solution of normality equation
\thetag{6.1}.\enddemo
\head
\S\,11. System of multidimensional type
and specifically two-dimensional systems.
\endhead
\rightheadtext{\S\,11. Specifically two-dimensional systems.}
    The solutions of normality equation considered in\S\,9
correspond to Newtonian dynamical systems with force fields
\thetag{4.8} and \thetag{4.31}. If $A(v)=k\cdot v^2$ these
systems are {\bf metrizable} (see formula \thetag{5.3} in \S\,5).
For other choice of $A(v)$ they aren't metrizable, but they belong
to more wide class of systems of {\bf multidimensional type} (see
formula \thetag{5.4} and formulas \thetag{5.16}, \thetag{5.17} in
\S\,5). Parameter $A$ for the force fields \thetag{5.4} is
determined by formula \thetag{5.5}:
$$
A=\frac{h(W)-\bigl<\nabla W,\,\bold v\bigr>}{W_v}.\hskip -3em
\tag11.1
$$
In variable $x$, $y$, $v$, $\theta$ formula \thetag{11.1} looks like
$$
A=\frac{h(W)-v\,(W_x\,\cos\theta+W_y\,\sin\theta)}{W_v}.\hskip -3em
\tag11.2
$$
Remember that here $h=h(W)$ is a function of one variable, while
$W=W(x,y,v)$ is a function of three variables.\par
     In multidimensional case $n\geqslant 3$ functions \thetag{11.1}
exhaust all solutions of reduced system of normality equations
(see theorem~5.1 in \S\,5 above and thesis \cite{22}). Two-dimensional
case is quite different in this point. Here there is only one
reduced normality equation. It has solutions with {\bf no
multidimensional analogs}. Let's construct an example of such
essentially two-dimensional solution using the relationship
\thetag{10.7} from \S\,10. For this purpose we choose $\varPhi(I_1,I_2)
=(\cos(I_1)-I_2)^2-1/2$. Then one can solve the equation \thetag{10.7}
in explicit form and find the function $b=b(v,\theta)$:
$$
b=\frac{v^2\,\sin 2\theta+2\,v\,u\,\cos\theta+v
\sqrt{v^2+4\,u\,v\,\sin\theta+2\,u^2}}
{4\,u\,\,v\sin\theta+2\,u^2-v^2\,\cos2\theta}.\hskip -3em
\tag11.3
$$
According to \thetag{10.3}, function $b(v,\theta)$ is logarithmic
derivative of the required function $A(v,\theta)$. Hence 
$A=e^{I(v,\theta)}$, where 
$$
\aligned
I&=\int\frac{v^2\,\sin 2\theta+2\,v\,u\,\cos\theta+v}
{4\,u\,\,v\sin\theta+2\,u^2-v^2\,\cos2\theta}\,d\theta\,+\\
\vspace{2ex}
&+\int\frac{\sqrt{v^2+4\,u\,v\,\sin\theta+2\,u^2}}
{4\,u\,\,v\sin\theta+2\,u^2-v^2\,\cos2\theta}\,d\theta.
\endaligned
\hskip -3em
\tag11.4
$$
First of two integrals \thetag{11.4} can be calculated explicitly in
elementary functions. Second one, upon the change of variable $\tau
=\sin\theta$, is brought to elliptic integrals (see \cite{Gra1}) and
expressed through elliptic functions (see \cite{Ahi1}). The dependence
on theta expressed by elliptic functions cannot be obtained by
formula \thetag{11.2}, no matter what function $W=W(x,y,v)$ is
chosen in it. The variable $\theta$ enters the formula \thetag{11.2}
in purely trigonometric form through $\sin\theta$ and $\cos\theta$.
This means that the solution of normality equation given by integrals
\thetag{11.4} is not among solutions of multidimensional type. It is
essentially two-dimensional solution.
\enddocument
\end